\documentclass[11pt]{article}
\usepackage{amsmath,amscd,amssymb}
\usepackage[english]{babel}
 \usepackage[applemac]{inputenc}

\newtheorem{Lem}{Lemma \thesection.}
\newtheorem{Th}[Lem]{Theorem \thesection.}
\newtheorem{Cor}[Lem]{Corollary \thesection.}
\newtheorem{Def}[Lem]{Definition \thesection.}
\newtheorem{Ex}[Lem]{Examples \thesection.}
\newtheorem{Prop}[Lem]{Proposition \thesection.}
\newtheorem{Rem}[Lem]{Remark \thesection.}
\newtheorem{Lem and Def}[Lem]{Lemma and Definition \thesection.}
\def\cal{\mathcal}
\def\bb{\mathbb} 
\def\a{\alpha }

\def\g{\gamma }
\def\d{\delta }
\def\D{\Delta }

\def\g{\gamma }
\def\G{\Gamma }
\def\k{\kappa }
\def\l{\lambda }

\def\o{\omega }
\def\p{\pi }
\def\P{\Pi }
\def\s{\sigma }
\def\S{\Sigma }
\def\t{\theta }
\def\T{\Theta }
\def\f{\varphi }
\def\dim{\rm dim\; }

\def\iff{ if and only if }

\def\Card{{\rm Card}}
\def\tr{\rm tr}

\newcommand{\tlowername}[2]%
{$\stackrel{\makebox[1pt]{#1}}%
{\begin{picture}(0,0)%
\put(0,0){\makebox(0,6)[t]{\makebox[1pt]{$#2$}}}%
\end{picture}}$}%

\newcommand{\AR}[1]%
{\begin{picture}(#1,0)%
\put(0,0){\vector(1,0){#1}}%
\end{picture}}%

\newcommand{\DOTAR}[1]%
{\NUMBEROFDOTS=#1%
\divide\NUMBEROFDOTS by 3%
\begin{picture}(#1,0)%
\multiput(0,0)(3,0){\NUMBEROFDOTS}{\circle*{1}}%
\put(#1,0){\vector(1,0){0}}%
\end{picture}}%

\newcommand{\MONO}[1]%
{\begin{picture}(#1,0)%
\put(0,0){\vector(1,0){#1}}%
\put(2,-2){\line(0,1){4}}%
\end{picture}}%

\newcommand{\EPI}[1]%
{\begin{picture}(#1,0)(-#1,0)%
\put(-#1,0){\vector(1,0){#1}}%
\put(-6,-2){\line(0,1){4}}%
\end{picture}}%

\newcommand{\BIMO}[1]%
{\begin{picture}(#1,0)(-#1,0)%
\put(-#1,0){\vector(1,0){#1}}%
\put(-6,-2){\line(0,1){4}}%
\put(-#1,-2){\hspace{2pt}\line(0,1){4}}%
\end{picture}}%

\newcommand{\BIAR}[1]%
{\begin{picture}(#1,4)%
\put(0,0){\vector(1,0){#1}}%
\put(0,4){\vector(1,0){#1}}%
\end{picture}}%

\newcommand{\EQL}[1]%
{\begin{picture}(#1,0)%
\put(0,1){\line(1,0){#1}}%
\put(0,-1){\line(1,0){#1}}%
\end{picture}}%

\newcommand{\ADJAR}[1]%
{\begin{picture}(#1,4)%
\put(0,0){\vector(1,0){#1}}%
\put(#1,4){\vector(-1,0){#1}}%
\end{picture}}%

\newcommand{\BKAR}[1]%
{\begin{picture}(#1,0)%
\put(#1,0){\vector(-1,0){#1}}%
\end{picture}}%

\newcommand{\BKDOTAR}[1]%
{\NUMBEROFDOTS=#1%
\divide\NUMBEROFDOTS by 3%
\begin{picture}(#1,0)%
\multiput(#1,0)(-3,0){\NUMBEROFDOTS}{\circle*{1}}%
\put(0,0){\vector(-1,0){0}}%
\end{picture}}%

\newcommand{\BKMONO}[1]%
{\begin{picture}(#1,0)(-#1,0)%
\put(0,0){\vector(-1,0){#1}}%
\put(-2,-2){\line(0,1){4}}%
\end{picture}}%

\newcommand{\BKEPI}[1]%
{\begin{picture}(#1,0)%
\put(#1,0){\vector(-1,0){#1}}%
\put(6,-2){\line(0,1){4}}%
\end{picture}}%

\newcommand{\BKBIMO}[1]%
{\begin{picture}(#1,0)%
\put(#1,0){\vector(-1,0){#1}}%
\put(6,-2){\line(0,1){4}}%
\put(#1,-2){\hspace{-2pt}\line(0,1){4}}%
\end{picture}}%

\newcommand{\BKBIAR}[1]%
{\begin{picture}(#1,4)%
\put(#1,0){\vector(-1,0){#1}}%
\put(#1,4){\vector(-1,0){#1}}%
\end{picture}}%

\newcommand{\BKADJAR}[1]%
{\begin{picture}(#1,4)%
\put(0,4){\vector(1,0){#1}}%
\put(#1,0){\vector(-1,0){#1}}%
\end{picture}}%

\newcommand{\lowername}[2]%
{$\stackrel{\makebox[1pt]{#1}}%
{\begin{picture}(0,0)%
\truex{600}%
\put(0,0){\makebox(0,\value{x})[t]{\makebox[1pt]{$#2$}}}%
\end{picture}}$}%

\newcommand{\hcase}[2]%
{\makebox[0pt]%
{\raisebox{-1pt}[0pt][0pt]{#1{#2}}}}%

\newcommand{\Hcase}[3]%
{\makebox[0pt]
{\raisebox{-1pt}[0pt][0pt]%
{$\stackrel{\makebox[0pt]{$\textstyle{#2}$}}{#1{#3}}$}}}%

\newcommand{\hcasE}[3]%
{\makebox[0pt]%
{\raisebox{-9pt}[0pt][0pt]%
{\lowername{#1{#3}}{#2}}}}%

\newcommand{\hbicase}[2]%
{\makebox[0pt]%
{\raisebox{-2.5pt}[0pt][0pt]{#1{#2}}}}%

\newcommand{\Hbicase}[4]%
{\makebox[0pt]
{\raisebox{-10.5pt}[0pt][0pt]%
{$\stackrel{\makebox[0pt]{$\textstyle{#2}$}}%
{\mbox{\lowername{#1{#4}}{#3}}}$}}}%

\newcommand{\EAR}[1]%
{\begin{picture}(#1,0)%
\put(0,0){\vector(1,0){#1}}%
\end{picture}}%

\newcommand{\EDOTAR}[1]%
{\truex{100}\truey{300}%
\NUMBEROFDOTS=#1%
\divide\NUMBEROFDOTS by \value{y}%
\begin{picture}(#1,0)%
\multiput(0,0)(\value{y},0){\NUMBEROFDOTS}%
{\circle*{\value{x}}}%
\put(#1,0){\vector(1,0){0}}%
\end{picture}}%

\newcommand{\EMONO}[1]%
{\begin{picture}(#1,0)%
\put(0,0){\vector(1,0){#1}}%
\truex{300}\truey{600}%
\put(\value{x},-\value{x}){\line(0,1){\value{y}}}%
\end{picture}}%

\newcommand{\EEPI}[1]%
{\begin{picture}(#1,0)(-#1,0)%
\put(-#1,0){\vector(1,0){#1}}%
\truex{300}\truey{600}\truez{800}%
\put(-\value{z},-\value{x}){\line(0,1){\value{y}}}%
\end{picture}}%

\newcommand{\EBIMO}[1]%
{\begin{picture}(#1,0)(-#1,0)%
\put(-#1,0){\vector(1,0){#1}}%
\truex{300}\truey{600}\truez{800}%
\put(-\value{z},-\value{x}){\line(0,1){\value{y}}}%
\put(-#1,-\value{x}){\hspace{3pt}\line(0,1){\value{y}}}%
\end{picture}}%

\newcommand{\EBIAR}[1]%
{\truex{400}%
\begin{picture}(#1,\value{x})%
\put(0,0){\vector(1,0){#1}}%
\put(0,\value{x}){\vector(1,0){#1}}%
\end{picture}}%

\newcommand{\EEQL}[1]%
{\begin{picture}(#1,0)%
\truex{200}%
\put(0,\value{x}){\line(1,0){#1}}%
\put(0,0){\line(1,0){#1}}%
\end{picture}}%

\newcommand{\EADJAR}[1]%
{\truex{400}%
\begin{picture}(#1,\value{x})%
\put(0,0){\vector(1,0){#1}}%
\put(#1,\value{x}){\vector(-1,0){#1}}%
\end{picture}}%

\newcommand{\ear}%
{\hspace{\SOURCE\unitlength}%
\hcase{\EAR}{\ARROWLENGTH}}%

\newcommand{\Ear}[1]%
{\hspace{\SOURCE\unitlength}%
\Hcase{\EAR}{#1}{\ARROWLENGTH}}%

\newcommand{\eaR}[1]%
{\hspace{\SOURCE\unitlength}%
\hcasE{\EAR}{#1}{\ARROWLENGTH}}%

\newcommand{\edotar}%
{\hspace{\SOURCE\unitlength}%
\hcase{\EDOTAR}{\ARROWLENGTH}}%

\newcommand{\Edotar}[1]%
{\hspace{\SOURCE\unitlength}%
\Hcase{\EDOTAR}{#1}{\ARROWLENGTH}}%

\newcommand{\edotaR}[1]%
{\hspace{\SOURCE\unitlength}%
\hcasE{\EDOTAR}{#1}{\ARROWLENGTH}}%

\newcommand{\emono}%
{\hspace{\SOURCE\unitlength}%
\hcase{\EMONO}{\ARROWLENGTH}}%

\newcommand{\Emono}[1]%
{\hspace{\SOURCE\unitlength}%
\Hcase{\EMONO}{#1}{\ARROWLENGTH}}%

\newcommand{\emonO}[1]%
{\hspace{\SOURCE\unitlength}%
\hcasE{\EMONO}{#1}{\ARROWLENGTH}}%

\newcommand{\eepi}%
{\hspace{\SOURCE\unitlength}%
\hcase{\EEPI}{\ARROWLENGTH}}%

\newcommand{\Eepi}[1]%
{\hspace{\SOURCE\unitlength}%
\Hcase{\EEPI}{#1}{\ARROWLENGTH}}%

\newcommand{\eepI}[1]%
{\hspace{\SOURCE\unitlength}%
\hcasE{\EEPI}{#1}{\ARROWLENGTH}}%

\newcommand{\ebimo}%
{\hspace{\SOURCE\unitlength}%
\hcase{\EBIMO}{\ARROWLENGTH}}%

\newcommand{\Ebimo}[1]%
{\hspace{\SOURCE\unitlength}%
\Hcase{\EBIMO}{#1}{\ARROWLENGTH}}%

\newcommand{\ebimO}[1]%
{\hspace{\SOURCE\unitlength}%
\hcasE{\EBIMO}{#1}{\ARROWLENGTH}}%

\newcommand{\eiso}%
{\hspace{\SOURCE\unitlength}%
\Hcase{\EAR}{\cong}{\ARROWLENGTH}}%

\newcommand{\Eiso}[1]%
{\hspace{\SOURCE\unitlength}%
\Hcase{\EAR}{\cong#1}{\ARROWLENGTH}}%

\newcommand{\eisO}[1]%
{\hspace{\SOURCE\unitlength}%
\hcasE{\EAR}{\cong#1}{\ARROWLENGTH}}%

\newcommand{\ebiar}%
{\hspace{\SOURCE\unitlength}%
\hbicase{\EBIAR}{\ARROWLENGTH}}%

\newcommand{\Ebiar}[2]%
{\hspace{\SOURCE\unitlength}%
\Hbicase{\EBIAR}{#1}{#2}{\ARROWLENGTH}}%

\newcommand{\eeql}%
{\hspace{\SOURCE\unitlength}%
\hbicase{\EEQL}{\ARROWLENGTH}}%

\newcommand{\eadjar}%
{\hspace{\SOURCE\unitlength}%
\hbicase{\EADJAR}{\ARROWLENGTH}}%

\newcommand{\Eadjar}[2]%
{\hspace{\SOURCE\unitlength}%
\Hbicase{\EADJAR}{#1}{#2}{\ARROWLENGTH}}%

\newcommand{\WAR}[1]%
{\begin{picture}(#1,0)%
\put(#1,0){\vector(-1,0){#1}}%
\end{picture}}%

\newcommand{\WDOTAR}[1]%
{\truex{100}\truey{300}%
\NUMBEROFDOTS=#1%
\divide\NUMBEROFDOTS by \value{y}%
\begin{picture}(#1,0)%
\multiput(#1,0)(-\value{y},0){\NUMBEROFDOTS}%
{\circle*{\value{x}}}%
\put(0,0){\vector(-1,0){0}}%
\end{picture}}%

\newcommand{\WMONO}[1]%
{\begin{picture}(#1,0)(-#1,0)%
\put(0,0){\vector(-1,0){#1}}%
\truex{300}\truey{600}%
\put(-\value{x},-\value{x}){\line(0,1){\value{y}}}%
\end{picture}}%

\newcommand{\WEPI}[1]%
{\begin{picture}(#1,0)%
\put(#1,0){\vector(-1,0){#1}}%
\truex{300}\truey{600}\truez{800}%
\put(\value{z},-\value{x}){\line(0,1){\value{y}}}%
\end{picture}}%

\newcommand{\WBIMO}[1]%
{\begin{picture}(#1,0)%
\put(#1,0){\vector(-1,0){#1}}%
\truex{300}\truey{600}\truez{800}%
\put(\value{z},-\value{x}){\line(0,1){\value{y}}}%
\put(#1,-\value{x}){\hspace{-3pt}\line(0,1){\value{y}}}%
\end{picture}}%

\newcommand{\WBIAR}[1]%
{\truex{400}%
\begin{picture}(#1,\value{x})%
\put(#1,0){\vector(-1,0){#1}}%
\put(#1,\value{x}){\vector(-1,0){#1}}%
\end{picture}}%

\newcommand{\WADJAR}[1]%
{\truex{400}%
\begin{picture}(#1,\value{x})%
\put(0,\value{x}){\vector(1,0){#1}}%
\put(#1,0){\vector(-1,0){#1}}%
\end{picture}}%

\newcommand{\war}%
{\hspace{\SOURCE\unitlength}%
\hcase{\WAR}{\ARROWLENGTH}}%

\newcommand{\War}[1]%
{\hspace{\SOURCE\unitlength}%
\Hcase{\WAR}{#1}{\ARROWLENGTH}}%

\newcommand{\waR}[1]%
{\hspace{\SOURCE\unitlength}%
\hcasE{\WAR}{#1}{\ARROWLENGTH}}%

\newcommand{\wdotar}%
{\hspace{\SOURCE\unitlength}%
\hcase{\WDOTAR}{\ARROWLENGTH}}%

\newcommand{\Wdotar}[1]%
{\hspace{\SOURCE\unitlength}%
\Hcase{\WDOTAR}{#1}{\ARROWLENGTH}}%

\newcommand{\wdotaR}[1]%
{\hspace{\SOURCE\unitlength}%
\hcasE{\WDOTAR}{#1}{\ARROWLENGTH}}%

\newcommand{\wmono}%
{\hspace{\SOURCE\unitlength}%
\hcase{\WMONO}{\ARROWLENGTH}}%

\newcommand{\Wmono}[1]%
{\hspace{\SOURCE\unitlength}%
\Hcase{\WMONO}{#1}{\ARROWLENGTH}}%

\newcommand{\wmonO}[1]%
{\hspace{\SOURCE\unitlength}%
\hcasE{\WMONO}{#1}{\ARROWLENGTH}}%

\newcommand{\wepi}%
{\hspace{\SOURCE\unitlength}%
\hcase{\WEPI}{\ARROWLENGTH}}%

\newcommand{\Wepi}[1]%
{\hspace{\SOURCE\unitlength}%
\Hcase{\WEPI}{#1}{\ARROWLENGTH}}%

\newcommand{\wepI}[1]%
{\hspace{\SOURCE\unitlength}%
\hcasE{\WEPI}{#1}{\ARROWLENGTH}}%

\newcommand{\wbimo}%
{\hspace{\SOURCE\unitlength}%
\hcase{\WBIMO}{\ARROWLENGTH}}%

\newcommand{\Wbimo}[1]%
{\hspace{\SOURCE\unitlength}%
\Hcase{\WBIMO}{#1}{\ARROWLENGTH}}%

\newcommand{\wbimO}[1]%
{\hspace{\SOURCE\unitlength}%
\hcasE{\WBIMO}{#1}{\ARROWLENGTH}}%

\newcommand{\wiso}%
{\hspace{\SOURCE\unitlength}%
\Hcase{\WAR}{\cong}{\ARROWLENGTH}}%

\newcommand{\Wiso}[1]%
{\hspace{\SOURCE\unitlength}%
\Hcase{\WAR}{#1}{\ARROWLENGTH}}%

\newcommand{\wisO}[1]%
{\hspace{\SOURCE\unitlength}%
\hcasE{\WAR}{#1}{\ARROWLENGTH}}%

\newcommand{\wbiar}%
{\hspace{\SOURCE\unitlength}%
\hbicase{\WBIAR}{\ARROWLENGTH}}%

\newcommand{\Wbiar}[2]%
{\hspace{\SOURCE\unitlength}%
\Hbicase{\WBIAR}{#1}{#2}{\ARROWLENGTH}}%

\newcommand{\weql}%
{\hspace{\SOURCE\unitlength}%
\hbicase{\EEQL}{\ARROWLENGTH}}%

\newcommand{\wadjar}%
{\hspace{\SOURCE\unitlength}%
\hbicase{\WADJAR}{\ARROWLENGTH}}%

\newcommand{\Wadjar}[2]%
{\hspace{\SOURCE\unitlength}%
\Hbicase{\WADJAR}{#1}{#2}{\ARROWLENGTH}}%

\newcommand{\Vcase}[3]{\makebox[0pt]%
{\makebox[0pt][r]{\raisebox{0pt}[0pt][0pt]{${#2}\hspace{2pt}$}}}#1{#3}}%

\newcommand{\vcasE}[3]{\makebox[0pt]%
{#1{#3}\makebox[0pt][l]{\raisebox{0pt}[0pt][0pt]{\hspace{2pt}$#2$}}}}%

\newcommand{\Vbicase}[4]{\makebox[0pt]%
{\makebox[0pt][r]{\raisebox{0pt}[0pt][0pt]{$#2$\hspace{4pt}}}#1{#4}%
\makebox[0pt][l]{\raisebox{0pt}[0pt][0pt]{\hspace{5pt}$#3$}}}}%

\newcommand{\SAR}[1]%
{\begin{picture}(0,0)%
\put(0,0){\makebox(0,0)%
{\begin{picture}(0,#1)%
\put(0,#1){\vector(0,-1){#1}}%
\end{picture}}}\end{picture}}%

\newcommand{\SDOTAR}[1]%
{\truex{100}\truey{300}%
\NUMBEROFDOTS=#1%
\divide\NUMBEROFDOTS by \value{y}%
\begin{picture}(0,0)%
\put(0,0){\makebox(0,0)%
{\begin{picture}(0,#1)%
\multiput(0,#1)(0,-\value{y}){\NUMBEROFDOTS}%
{\circle*{\value{x}}}%
\put(0,0){\vector(0,-1){0}}%
\end{picture}}}\end{picture}}%

\newcommand{\SMONO}[1]%
{\begin{picture}(0,0)%
\put(0,0){\makebox(0,0)%
{\begin{picture}(0,#1)%
\put(0,#1){\vector(0,-1){#1}}%
\truex{300}\truey{600}%
\put(0,#1){\begin{picture}(0,0)%
\put(-\value{x},-\value{x}){\line(1,0){\value{y}}}\end{picture}}%
\end{picture}}}\end{picture}}%

\newcommand{\SEPI}[1]%
{\begin{picture}(0,0)%
\put(0,0){\makebox(0,0)%
{\begin{picture}(0,#1)%
\put(0,#1){\vector(0,-1){#1}}%
\truex{300}\truey{600}\truez{800}%
\put(-\value{x},\value{z}){\line(1,0){\value{y}}}%
\end{picture}}}\end{picture}}%

\newcommand{\SBIMO}[1]%
{\begin{picture}(0,0)%
\put(0,0){\makebox(0,0)%
{\begin{picture}(0,#1)%
\put(0,#1){\vector(0,-1){#1}}%
\truex{300}\truey{600}\truez{800}%
\put(0,#1){\begin{picture}(0,0)%
\put(-\value{x},-\value{x}){\line(1,0){\value{y}}}\end{picture}}%
\put(-\value{x},\value{z}){\line(1,0){\value{y}}}%
\end{picture}}}\end{picture}}%

\newcommand{\SBIAR}[1]%
{\begin{picture}(0,0)%
\truex{200}%
\put(0,0){\makebox(0,0)%
{\begin{picture}(0,#1)\put(-\value{x},#1){\vector(0,-1){#1}}%
\put(\value{x},#1){\vector(0,-1){#1}}%
\end{picture}}}\end{picture}}%

\newcommand{\SEQL}[1]%
{\begin{picture}(0,0)%
\truex{100}%
\put(0,0){\makebox(0,0)%
{\begin{picture}(0,#1)\put(-\value{x},#1){\line(0,-1){#1}}%
\put(\value{x},#1){\line(0,-1){#1}}%
\end{picture}}}\end{picture}}%

\newcommand{\Sarv}[2]{\Vcase{\SAR}{#1}{#200}}%

\newcommand{\Sar}[1]{\Sarv{#1}{50}}%

\newcommand{\Sisov}[2]%
{\Vbicase{\SAR}{#1\hspace{-2pt}}{\hspace{-2pt}\cong}{#200}}%

\newcommand{\NAR}[1]%
{\begin{picture}(0,0)%
\put(0,0){\makebox(0,0)%
{\begin{picture}(0,#1)\put(0,0){\vector(0,1){#1}}%
\end{picture}}}\end{picture}}%

\newcommand{\NDOTAR}[1]%
{\truex{100}\truey{300}%
\NUMBEROFDOTS=#1%
\divide\NUMBEROFDOTS by \value{y}%
\begin{picture}(0,0)%
\put(0,0){\makebox(0,0)%
{\begin{picture}(0,#1)%
\multiput(0,0)(0,\value{y}){\NUMBEROFDOTS}%
{\circle*{\value{x}}}%
\put(0,#1){\vector(0,1){0}}%
\end{picture}}}\end{picture}}%

\newcommand{\NMONO}[1]%
{\begin{picture}(0,0)%
\put(0,0){\makebox(0,0)%
{\begin{picture}(0,#1)%
\put(0,0){\vector(0,1){#1}}%
\truex{300}\truey{600}%
\put(-\value{x},\value{x}){\line(1,0){\value{y}}}%
\end{picture}}}%
\end{picture}}%

\newcommand{\NEPI}[1]%
{\begin{picture}(0,0)%
\put(0,0){\makebox(0,0)%
{\begin{picture}(0,#1)%
\put(0,0){\vector(0,1){#1}}%
\truex{300}\truey{600}\truez{800}%
\put(0,#1){\begin{picture}(0,0)%
\put(-\value{x},-\value{z}){\line(1,0){\value{y}}}\end{picture}}%
\end{picture}}}\end{picture}}%

\newcommand{\NBIMO}[1]%
{\begin{picture}(0,0)%
\put(0,0){\makebox(0,0)%
{\begin{picture}(0,#1)%
\put(0,0){\vector(0,1){#1}}%
\truex{300}\truey{600}\truez{800}%
\put(-\value{x},\value{x}){\line(1,0){\value{y}}}%
\put(0,#1){\begin{picture}(0,0)%
\put(-\value{x},-\value{z}){\line(1,0){\value{y}}}\end{picture}}%
\end{picture}}}\end{picture}}%

\newcommand{\NBIAR}[1]%
{\begin{picture}(0,0)%
\truex{200}%
\put(0,0){\makebox(0,0)%
{\begin{picture}(0,#1)\put(-\value{x},0){\vector(0,1){#1}}%
\put(\value{x},0){\vector(0,1){#1}}%
\end{picture}}}\end{picture}}%

\newcommand{\naRv}[2]{\vcasE{\NAR}{#1}{#200}}%

\newcommand{\naR}[1]{\naRv{#1}{50}}%

\newcommand{\Nisov}[2]%
{\Vbicase{\NAR}{#1\hspace{-2pt}}{\hspace{-2pt}\cong}{#200}}%

\newcommand{\fdcase}[3]{\begin{picture}(0,0)%
\put(0,-150){#1}%
\truex{200}\truey{600}\truez{600}%
\put(-\value{x},-\value{x}){\makebox(0,\value{z})[r]{${#2}$}}%
\put(\value{x},-\value{y}){\makebox(0,\value{z})[l]{${#3}$}}%
\end{picture}}%

\newcommand{\NEAR}{\begin{picture}(0,0)%
\put(-2900,-2900){\vector(1,1){5800}}%
\end{picture}}%

\newcommand{\NEDOTAR}%
{\truex{100}\truey{212}%
\NUMBEROFDOTS=5800%
\divide\NUMBEROFDOTS by \value{y}%
\begin{picture}(0,0)%
\multiput(-2900,-2900)(\value{y},\value{y}){\NUMBEROFDOTS}%
{\circle*{\value{x}}}%
\put(2900,2900){\vector(1,1){0}}%
\end{picture}}%

\newcommand{\Near}[1]{\fdcase{\NEAR}{#1}{}}%

\newcommand{\SWAR}{\begin{picture}(0,0)%
\put(2900,2900){\vector(-1,-1){5800}}%
\end{picture}}%

\newcommand{\SWDOTAR}%
{\truex{100}\truey{212}%
\NUMBEROFDOTS=5800%
\divide\NUMBEROFDOTS by \value{y}%
\begin{picture}(0,0)%
\multiput(2900,2900)(-\value{y},-\value{y}){\NUMBEROFDOTS}%
{\circle*{\value{x}}}%
\put(-2900,-2900){\vector(-1,-1){0}}%
\end{picture}}%

\newcommand{\swaR}[1]{\fdcase{\SWAR}{}{#1}}%

\newcommand{\sdcase}[3]{\begin{picture}(0,0)%
\put(0,-150){#1}%
\truex{100}\truez{600}%
\put(\value{x},\value{x}){\makebox(0,\value{z})[l]{${#2}$}}%
\truex{300}\truey{800}%
\put(-\value{x},-\value{y}){\makebox(0,\value{z})[r]{${#3}$}}%
\end{picture}}%

\newcommand{\SEAR}{\begin{picture}(0,0)%
\put(-2900,2900){\vector(1,-1){5800}}%
\end{picture}}%

\newcommand{\SEDOTAR}%
{\truex{100}\truey{212}%
\NUMBEROFDOTS=5800%
\divide\NUMBEROFDOTS by \value{y}%
\begin{picture}(0,0)%
\multiput(-2900,2900)(\value{y},-\value{y}){\NUMBEROFDOTS}%
{\circle*{\value{x}}}%
\put(2900,-2900){\vector(1,-1){0}}%
\end{picture}}%

\newcommand{\Sear}[1]{\sdcase{\SEAR}{#1}{}}%

\newcommand{\seaR}[1]{\sdcase{\SEAR}{}{#1}}%

\newcommand{\NWDOTAR}%
{\truex{100}\truey{212}%
\NUMBEROFDOTS=5800%
\divide\NUMBEROFDOTS by \value{y}%
\begin{picture}(0,0)%
\multiput(2900,-2900)(-\value{y},\value{y}){\NUMBEROFDOTS}%
{\circle*{\value{x}}}%
\put(-2900,2900){\vector(-1,1){0}}%
\end{picture}}%

\newcommand{\ENEAR}[2]%
{\makebox[0pt]{\begin{picture}(0,0)%
\put(0,-150){\makebox(0,0){\begin{picture}(0,0)%
\put(-6600,-3300){\vector(2,1){13200}}%
\truex{200}\truey{800}\truez{600}%
\put(-\value{x},\value{x}){\makebox(0,\value{z})[r]{${#1}$}}%
\put(\value{x},-\value{y}){\makebox(0,\value{z})[l]{${#2}$}}%
\end{picture}}}\end{picture}}}%

\newcommand{\ESEAR}[2]%
{\makebox[0pt]{\begin{picture}(0,0)%
\put(0,-150){\makebox(0,0){\begin{picture}(0,0)%
\put(-6600,3300){\vector(2,-1){13200}}%
\truex{200}\truey{800}\truez{600}%
\put(\value{x},\value{x}){\makebox(0,\value{z})[l]{${#1}$}}%
\put(-\value{x},-\value{y}){\makebox(0,\value{z})[r]{${#2}$}}%
\end{picture}}}\end{picture}}}%

\newcommand{\WNWAR}[2]%
{\makebox[0pt]{\begin{picture}(0,0)%
\put(0,-150){\makebox(0,0){\begin{picture}(0,0)%
\put(6600,-3300){\vector(-2,1){13200}}%
\truex{200}\truey{800}\truez{600}%
\put(\value{x},\value{x}){\makebox(0,\value{z})[l]{${#1}$}}%
\put(-\value{x},-\value{y}){\makebox(0,\value{z})[r]{${#2}$}}%
\end{picture}}}\end{picture}}}%

\newcommand{\WSWAR}[2]%
{\makebox[0pt]{\begin{picture}(0,0)%
\put(0,-150){\makebox(0,0){\begin{picture}(0,0)%
\put(6600,3300){\vector(-2,-1){13200}}%
\truex{200}\truey{800}\truez{600}%
\put(-\value{x},\value{x}){\makebox(0,\value{z})[r]{${#1}$}}%
\put(\value{x},-\value{y}){\makebox(0,\value{z})[l]{${#2}$}}%
\end{picture}}}\end{picture}}}%

\newcommand{\NNEAR}[2]%
{\raisebox{-1pt}[0pt][0pt]{\begin{picture}(0,0)%
\put(0,0){\makebox(0,0){\begin{picture}(0,0)%
\put(-3300,-6600){\vector(1,2){6600}}%
\truex{100}\truez{600}%
\put(-\value{x},\value{x}){\makebox(0,\value{z})[r]{${#1}$}}%
\put(\value{x},-\value{z}){\makebox(0,\value{z})[l]{${#2}$}}%
\end{picture}}}\end{picture}}}%

\newcommand{\SSWAR}[2]%
{\raisebox{-1pt}[0pt][0pt]{\begin{picture}(0,0)%
\put(0,0){\makebox(0,0){\begin{picture}(0,0)%
\put(3300,6600){\vector(-1,-2){6600}}%
\truex{100}\truez{600}%
\put(-\value{x},\value{x}){\makebox(0,\value{z})[r]{${#1}$}}%
\put(\value{x},-\value{z}){\makebox(0,\value{z})[l]{${#2}$}}%
\end{picture}}}\end{picture}}}%

\newcommand{\SSEAR}[2]%
{\raisebox{-1pt}[0pt][0pt]{\begin{picture}(0,0)%
\put(0,0){\makebox(0,0){\begin{picture}(0,0)%
\put(-3300,6600){\vector(1,-2){6600}}%
\truex{200}\truez{600}%
\put(\value{x},\value{x}){\makebox(0,\value{z})[l]{${#1}$}}%
\put(-\value{x},-\value{z}){\makebox(0,\value{z})[r]{${#2}$}}%
\end{picture}}}\end{picture}}}%

\newcommand{\NNWAR}[2]%
{\raisebox{-1pt}[0pt][0pt]{\begin{picture}(0,0)%
\put(0,0){\makebox(0,0){\begin{picture}(0,0)%
\put(3300,-6600){\vector(-1,2){6600}}%
\truex{200}\truez{600}%
\put(\value{x},\value{x}){\makebox(0,\value{z})[l]{${#1}$}}%
\put(-\value{x},-\value{z}){\makebox(0,\value{z})[r]{${#2}$}}%
\end{picture}}}\end{picture}}}%

\newcommand{\Necurve}[2]%
{\begin{picture}(0,0)%
\truex{1300}\truey{2000}\truez{200}%
\put(0,\value{x}){\oval(#200,\value{y})[t]}%
\put(0,\value{x}){\makebox(0,0){\begin{picture}(#200,0)%
\put(#200,0){\vector(0,-1){\value{z}}}%
\put(0,0){\line(0,-1){\value{z}}}\end{picture}}}%
\truex{2500}%
\put(0,\value{x}){\makebox(0,0)[b]{${#1}$}}%
\end{picture}}%

\newcommand{\Nwcurve}[2]%
{\begin{picture}(0,0)%
\truex{1300}\truey{2000}\truez{200}%
\put(0,\value{x}){\oval(#200,\value{y})[t]}%
\put(0,\value{x}){\makebox(0,0){\begin{picture}(#200,0)%
\put(#200,0){\line(0,-1){\value{z}}}%
\put(0,0){\vector(0,-1){\value{z}}}\end{picture}}}%
\truex{2500}%
\put(0,\value{x}){\makebox(0,0)[b]{${#1}$}}%
\end{picture}}%

\newcommand{\Securve}[2]%
{\begin{picture}(0,0)%
\truex{1300}\truey{2000}\truez{200}%
\put(0,-\value{x}){\oval(#200,\value{y})[b]}%
\put(0,-\value{x}){\makebox(0,0){\begin{picture}(#200,0)%
\put(#200,0){\vector(0,1){\value{z}}}%
\put(0,0){\line(0,1){\value{z}}}\end{picture}}}%
\truex{2500}%
\put(0,-\value{x}){\makebox(0,0)[t]{${#1}$}}%
\end{picture}}%

\newcommand{\Swcurve}[2]%
{\begin{picture}(0,0)%
\truex{1300}\truey{2000}\truez{200}%
\put(0,-\value{x}){\oval(#200,\value{y})[b]}%
\put(0,-\value{x}){\makebox(0,0){\begin{picture}(#200,0)%
\put(#200,0){\line(0,1){\value{z}}}%
\put(0,0){\vector(0,1){\value{z}}}\end{picture}}}%
\truex{2500}%
\put(0,-\value{x}){\makebox(0,0)[t]{${#1}$}}%
\end{picture}}%

\newcommand{\Escurve}[2]%
{\begin{picture}(0,0)%
\truex{1400}\truey{2000}\truez{200}%
\put(\value{x},0){\oval(\value{y},#200)[r]}%
\put(\value{x},0){\makebox(0,0){\begin{picture}(0,#200)%
\put(0,0){\vector(-1,0){\value{z}}}%
\put(0,#200){\line(-1,0){\value{z}}}\end{picture}}}%
\truex{2500}%
\put(\value{x},0){\makebox(0,0)[l]{${#1}$}}%
\end{picture}}%

\newcommand{\Encurve}[2]%
{\begin{picture}(0,0)%
\truex{1400}\truey{2000}\truez{200}%
\put(\value{x},0){\oval(\value{y},#200)[r]}%
\put(\value{x},0){\makebox(0,0){\begin{picture}(0,#200)%
\put(0,0){\line(-1,0){\value{z}}}%
\put(0,#200){\vector(-1,0){\value{z}}}\end{picture}}}%
\truex{2500}%
\put(\value{x},0){\makebox(0,0)[l]{${#1}$}}%
\end{picture}}%

\newcommand{\Wscurve}[2]%
{\begin{picture}(0,0)%
\truex{1300}\truey{2000}\truez{200}%
\put(-\value{x},0){\oval(\value{y},#200)[l]}%
\put(-\value{x},0){\makebox(0,0){\begin{picture}(0,#200)%
\put(0,0){\vector(1,0){\value{z}}}%
\put(0,#200){\line(1,0){\value{z}}}\end{picture}}}%
\truex{2400}%
\put(-\value{x},0){\makebox(0,0)[r]{${#1}$}}%
\end{picture}}%

\newcommand{\Wncurve}[2]%
{\begin{picture}(0,0)%
\truex{1300}\truey{2000}\truez{200}%
\put(-\value{x},0){\oval(\value{y},#200)[l]}%
\put(-\value{x},0){\makebox(0,0){\begin{picture}(0,#200)%
\put(0,0){\line(1,0){\value{z}}}%
\put(0,#200){\vector(1,0){\value{z}}}\end{picture}}}%
\truex{2400}%
\put(-\value{x},0){\makebox(0,0)[r]{${#1}$}}%
\end{picture}}%

\newcount\SCALE%

\newcount\NUMBER%

\newcount\LINE%

\newcount\COLUMN%

\newcount\WIDTH%

\newcount\SOURCE%

\newcount\ARROW%

\newcount\TARGET%

\newcount\ARROWLENGTH%

\newcount\NUMBEROFDOTS%

\newcounter{x}%

\newcounter{y}%

\newcounter{z}%

\newcounter{horizontal}%

\newcounter{vertical}%

\newskip\itemlength%

\newskip\firstitem%

\newskip\seconditem%

\newcommand{\printarrow}{}%

\newcommand{\truex}[1]{%
\NUMBER=#1%
\multiply\NUMBER by 100%
\divide\NUMBER by \SCALE%
\setcounter{x}{\NUMBER}}%

\newcommand{\truey}[1]{%
\NUMBER=#1%
\multiply\NUMBER by 100%
\divide\NUMBER by \SCALE%
\setcounter{y}{\NUMBER}}%

\newcommand{\truez}[1]{%
\NUMBER=#1%
\multiply\NUMBER by 100%
\divide\NUMBER by \SCALE%
\setcounter{z}{\NUMBER}}%

\newcommand{\changecounters}[1]{%
\SOURCE=\ARROW%
\ARROW=\TARGET%
\settowidth{\itemlength}{#1}%
\ifdim \itemlength > 2800\unitlength%
\addtolength{\itemlength}{-2800\unitlength}%
\TARGET=\itemlength%
\divide\TARGET by 1310%
\multiply\TARGET by 100%
\divide\TARGET by \SCALE%
\else%
\TARGET=0%
\fi%
\ARROWLENGTH=5000%
\advance\ARROWLENGTH by -\SOURCE%
\advance\ARROWLENGTH by -\TARGET%
\advance\SOURCE by -\TARGET}%

\newcommand{\initialize}[1]{%
\LINE=0%
\COLUMN=0%
\WIDTH=0%
\ARROW=0%
\TARGET=0%
\changecounters{#1}%
\renewcommand{\printarrow}{#1}%
\begin{center}%
\vspace{10pt}%
\begin{picture}(0,0)}%

\newcommand{\DIAGV}[2]{%
\SCALE=#1%
\setlength{\unitlength}{655sp}%
\multiply\unitlength by \SCALE%
\divide\unitlength by 100%
\initialize{\mbox{$#2$}}}%

\newcommand{\n}[1]{%
\changecounters{\mbox{$#1$}}%
\put(\COLUMN,\LINE){\makebox(0,0){\printarrow}}%
\thinlines%
\renewcommand{\printarrow}{\mbox{$#1$}}%
\advance\COLUMN by 4000}%

\newcommand{\nn}[1]{%
\put(\COLUMN,\LINE){\makebox(0,0){\printarrow}}%
\thinlines%
\ifnum \WIDTH < \COLUMN%
\WIDTH=\COLUMN%
\else%
\fi%
\advance\LINE by -4000%
\COLUMN=0%
\ARROW=0%
\TARGET=0%
\changecounters{\mbox{$#1$}}%
\renewcommand{\printarrow}{\mbox{$#1$}}}%

\newcommand{\conclude}{%
\put(\COLUMN,\LINE){\makebox(0,0){\printarrow}}%
\thinlines%
\ifnum \WIDTH < \COLUMN%
\WIDTH=\COLUMN%
\else%
\fi%
\setcounter{horizontal}{\WIDTH}%
\setcounter{vertical}{-\LINE}%
\end{picture}}%

\newcommand{\diag}{%
\conclude%
\raisebox{0pt}[0pt][\value{vertical}\unitlength]{}%
\hspace*{\value{horizontal}\unitlength}%
\vspace{10pt}%
\end{center}%
\setlength{\unitlength}{1pt}}%

\newcommand{\diagv}[3]{%
\conclude%
\NUMBER=#1%
\rule{0pt}{\NUMBER pt}%
\hspace*{-#2pt}%
\raisebox{0pt}[0pt][\value{vertical}\unitlength]{}%
\hspace*{\value{horizontal}\unitlength}
\NUMBER=#3%
\advance\NUMBER by 10%
\vspace*{\NUMBER pt}%
\end{center}%
\setlength{\unitlength}{1pt}}%

\newcommand{\N}[1]%
{\raisebox{0pt}[7pt][0pt]{$#1$}}%

\newcommand{\crosslength}[2]{%
\settowidth{\firstitem}{#1}%
\settowidth{\seconditem}{#2}%
\ifdim\firstitem < \seconditem%
\itemlength=\seconditem%
\else%
\itemlength=\firstitem%
\fi%
\divide\itemlength by 2%
\hspace{\itemlength}}%

\begin{document}

 \title {On surfaces of class VII$_0^+$ with numerically anticanonical divisor}
\author {Georges Dloussky}
\date{ }
\maketitle
\begin{abstract}We consider minimal compact complex surfaces $S$ with Betti numbers $b_1=1$ and $n=b_2>0$. A theorem of Donaldson gives $n$ exceptional line bundles. We prove that if in a deformation, these line bundles have sections, $S$ is a degeneration of blown-up Hopf surfaces. Besides, if there exists an integer $m\ge 1$ and a flat line bundle $F$ such that $H^0(S,-mK\otimes F)\neq 0$, then $S$ contains a Global Spherical Shell. We apply this last result to complete classification of bihermitian surfaces.
\end{abstract}
\tableofcontents
\section{Introduction}
A  minimal compact complex surface $S$ is said to be  of the class VII$_0$  of Kodaira if the first Betti number satisfies $b_1(S)=1$. A surface $S$ is of class VII$_0^+$ if moreover 
$n:=b_2(S)>0$; these surfaces admit no nonconstant meromorphic functions. 
The major problem in classification of non-k\"ahlerian surfaces is to achieve the classification of surfaces $S$ of class VII$_0^+$. All known surfaces of this class contain Global Spherical Shells (GSS), i.e. admit a biholomorphic map $\f:U\to V$ from a neighbourhood $U\subset \bb C^2\setminus\{0\}$ of the sphere $S^3=\partial B^2$  onto an open set $V$ such that $\S=\f(S^3)$ does not disconnect $S$. Are there other surfaces~? 
In first section we investigate the general situation: A theorem of Donaldson \cite{DON} gives a $\bb Z$-base $(E_i)$ of $H^2(S,\bb Z)$, such that $E_iE_j=-\d_{ij}$. These cohomology classes can be represented by line bundles $L_i$ such that $K_SL_i=L_i^2=-1$. Indeed, these line bundles generalize exceptional curves of the first kind, and since $S$ is minimal,  they have no section. Over the versal deformation $\cal S\to B$ of $S$ these line bundles form families $\cal L_i$.   We propose the following conjecture which can be easily checked for surfaces with GSS:\\
{\bf Conjecture 1}: Let $S$ be a surface in class VII$_0^+$ and $\cal S\to B$ be the versal deformation of $S\simeq S_0$ over the ball of dimension $h^1(S,\T)$. Then there exists $u\in \D$, $u\neq 0$, and flat line bundles $F_i$ such that $H^0(S_u, L_{i,u}\otimes F_i)\neq 0$ for $i=0,\ldots,n-1$.\\
 The main result of section 1 is (see theorem (\ref{Th1})),
 \begin{Th} Let $S$ be a surface in class VII$_0^+$ and $\cal S\to B$ its versal deformation. If there exists $u\in B$ and flat line bundles $F_i\in H^1(S,\bb C^\star)$ such that $H^0(S_u, L_{i,u}\otimes F_i)\neq 0$ for $i=0,\ldots,n-1$, then there is a non empty Zariski open set $U\subset B$ such that for all $u\in U$, $S_u$ is a blown-up Hopf surface. In particular, $S$ is a degeneration of blown-up Hopf surfaces.
 \end{Th}
If a surface is a degeneration of blown-up Hopf surfaces, the fundamental group of a fiber is isomorphic to $\bb Z\times \bb Z_l$, hence taking a finite covering, once obtains a surface obtained by degeneration of blown-up primary Hopf surfaces. Notice that a finite quotient of a surface of class VII$_0^+$ containing a GSS still contains a GSS \cite{D3}.\\ \\
{\bf Conjecture 2}: Let $S$ be a surface of class VII$_0^+$. If $S$ is a degeneration of blown-up primary Hopf surfaces, then $S$ contains a cycle of rational curves.\\ \\
A surface admitting a numerically anticanonical (NAC) divisor (see (\ref{DefNAC})), contains a cycle of rational curves. In section 2, we shall prove
\begin{Th} Let $S$ be a surface of class VII$_0^+$. If $S$ admits a NAC divisor, then $S$ contains a GSS.
\end{Th}
It is a weak version of\\
{\bf Conjecture 3 (Nakamura \cite{N2})}. Let $S$ be a surface of class VII$_0^+$. If $S$ contains a cycle of rational curves, $S$ contains a GSS.\\
The proof is based on the fact that in $H^2(S,\bb Z)$, a curve is equivalent to  a class of the form $L_i-\sum_{j\in I}L_j$, with $I\neq \emptyset$. Intuitively $L_i$ represents an exceptional curve of the first kind and $C$ is then equivalent to an exceptional curve of the first kind blown-up several times (Card$(I)$ times). It explains why curves have self-intersection $\le -2$. We recover a characterization of Inoue-Hirzebruch surfaces by Oeljeklaus, Toma \& Zaffran \cite{OTZ}:
\begin{Th} Let $S$ be a surface of class VII$_0$ with $b_2(S)>0$. Then $S$ is a Inoue-Hirzebruch surface if and only if  there exists  two flat line bundles $F_1$, $F_2$, two twisted vector fields $\t_1\in H^0(S,\T\otimes F_1)$, $\t_2\in H^0(S,\T\otimes F_2)$,   such that $\t_1\wedge\t_2(p)\neq 0$ at at least one point $p\in S$.
\end{Th}
In section 3 we apply  results of section 2 to complete the classification of bihermitian 4-manifolds $M$ (see \cite{A},\cite{AGG} \cite{PON}), when $b_1(M)=1$ and $b_2(M)>0$: A bihermitian surface is a riemannian oriented connected 4-manifold $(M,g)$ endowed with two integrable almost complex structures $J_1$, $J_2$ inducing the same orientation, orthogonal with respect to $g$ and independent i.e. $J_1(x)\neq \pm J_2(x)$ for at least one point $x\in M$. This structure depends only on the conformal class $c$ of $g$. A bihermitian surface is strongly bihermitian if $J_1(x)\neq \pm J_2(x)$ for every point $x\in M$. The key observation is that under these assumptions, $(M,J_i)$, $i=1,2$ admit a numerically anticanonical divisor.
\begin{Th} Let $(M,c,J_1,J_2)$ be a compact bihermitian surface with odd first Betti number. \\
1) If $(M,c,J_1,J_2)$ is strongly bihermitian (i.e $\cal D=\emptyset$), then the complex surfaces $(M,J_i)$ are minimal and either a Hopf surface covered by a primary one associated to a contraction $F:(\bb C^2,0)\to (\bb C^2,0)$ of the form
$$F(z_1,z_2)=(\a z_1+sz_2^m, a\a^{-1}z_2), $$
$$ {\rm with}\quad a, s\in\bb C, \ 0<|\a|^2\le a< |\a|<1, \ (a^m-\a^{m+1})s=0,$$
or else $(M,J_i)$ are Inoue surfaces $S^+_{N,p,q,r;t}$, $S^-_{N,p,q,r}$.\\
2) If $(M,c,J_1,J_2)$ is not strongly bihermitian, then $\cal D$ has at most two connected components, $(M,J_i)$, $i=1,2$, contain GSS and the minimal models $S_i$ of $(M,J_i)$ are
\begin{itemize}
\item Surfaces with GSS of intermediate type if $\cal D$ has one connected component
\item Hopf surfaces of special type (see {\rm \cite{PON} 2.2}), Inoue (parabolic) surfaces or Inoue-Hirzebruch surfaces if $\cal D$ has two connected components.
\end{itemize}
Moreover, the blown-up points belong to the NAC divisors.
\end{Th}
If moreover the metric $g$ is anti-self-dual (ASD), we obtain

\begin{Cor}  Let $(M,c,J_1,J_2)$ be a compact  ASD bihermitian surface with odd first Betti number. Then the minimal models of the complex surfaces $(M,J_i)$, $i=1,2$, are 
\begin{itemize}
\item Hopf surfaces of special type (see \cite{PON} 2.2),
\item (parabolic) Inoue surfaces or
\item even  Inoue-Hirzebruch surfaces.
\end{itemize}
Moreover, the blown-up points belong to the NAC divisors.
\end{Cor}

\begin{Rem} {\rm Throughout the paper we shall use the following terminology:\\
1) A surface for which exists a nontrivial divisor $D$ such that $D^2=0$ will be called a {\bf Enoki surface}, they are obtained by holomorphic compactification of an affine line bundle over an elliptic curve by a cycle $D$ of rational curves \cite{E}; otherwise they are associated to contracting holomorphic germs of maps 
$$F(z_1,z_2)=(t^nz_1z_2^n+\sum_{i=0}^{n-1}a_it^{i+1}z_2^{i+1},tz_2)$$
\cite{D1} \cite{DK}. A Enoki surface with an elliptic curve will be called briefly a {\bf Inoue surface (= parabolic Inoue surface)}: they are obtained by holomorphic  compactification of a line bundle over an elliptic curve by a cycle of rational curves $D$; otherwise they are associated to the contracting germs of maps
$$F(z_1,z_2)=(t^nz_1z_2^n,tz_2).$$
For all these surfaces, the sum of opposite self-intersections of the $n=b_2(S)$ rational curves $D_0,\ldots,D_{n-1}$ is $\s_n(S):=-\sum_{i=0}^{n-1}D_i^2=2n$.\\
2) A surface with $2n<\s_n(S)<3n$ will be called an {\bf intermediate surface};\\
3) A surface with $\s_n(S)=3n$ a called a {\bf Inoue-Hirzebruch (IH) surface} (\cite{I} or \cite{D2} for a construction by contracting germs of mappings). An {\bf even (=hyperbolic)} (resp. {\bf odd (=half)}) {\bf IH surfaces} has two (resp. one) cycle of rational curves.}
\end{Rem}

{\bf We shall assume throughout the article that $S$ admits no nonconstant meromorphic functions}.\\

{\bf Acknowledgments}: I am grateful for helpful discussions with Akira Fujiki, Massimiliano Pontecorvo,  and Vestislav Apostolov  during the preparation of the third part of this article. There is a gap in the proof of Corollary 2, p 425 in \cite{AGG}, theorems 4.1 and 5.2 in \cite{PON}. The statements are too restrictive and omit Inoue-Hirzebruch surfaces. 

\section{Exceptional line bundles and degeneration of blown-up Hopf surfaces}
Let $S$ be a surface in class VII$_0^+$ with $n=b_2(S)$. Since $S$ is not algebraic, $A^2\leq 0$  for every 
divisor $A$. By adjunction formula it is easy to deduce that  for every irreducible curve $C$, $KC\geq 0$, 
 $C^2\leq -2$ if $C$ is a regular rational
curve and  $C^2\leq 0$ if $C$ is a rational curve with a double 
point or an elliptic curve. It is well known that $S$ contains at most $n$ rational curves, and at most one elliptic curve.  By Hirzebruch index theorem, $b^-=b_2(S)$ whence the intersection form $Q:H^2(S,\bb Z)/Torsion \to \bb Z$ is negative definite.

\subsection{Exceptional line bundles}
An irreducible curve which satisfies the two conditions $KC=C^2=-1$ is an exceptional curve of the first kind; we generalize the notion to line bundles.
\begin{Def} 1) A line bundle $L\in H^1(S,\cal O^\star)$ is called an {\bf exceptional line bundle} (of the first kind) if $KL=L^2=-1$. \\
2) An effective divisor $E$ is called an exceptional divisor (of the first) kind if $E$ is the inverse image $\P^\star C$  of an exceptional curve of the first kind $C$ by a finite number of blowing-ups $\P$.
Equivalentlly it is an effective reduced divisor which may be blown-down onto a regular point.
\end{Def}
Using the fact that for a blowing-up $\P:S\to S'$, $K_S=\P^\star K_{S'}+C$ and the projection formula (\cite{BPV} p11), it is easy to check that if $E$ is an exceptional divisor, then $[E]$ is an exceptional line bundle. Moreover, the inverse image $\P^\star L$ of an exceptional line bundle $L$ by a finite sequence of blowing-ups is still an exceptional line bundle.\\

The following theorem has been proved by I. Nakamura when $S$ contains a cycle of rational curves \cite{N2} (1.7). It should be noticed that any surface with $b_1(S)=1$ and $b_2(S)=0$ is minimal and $H^1(S,\bb C^\star)\simeq H^1(S,\cal O^\star)$ \cite{KO} II, p699.
\begin{Th}\label{Zbase} Let $S$ be a  (not necessarily minimal) compact complex surface such that $b_1(S)=1$, with second Betti number $n=b_2(S)>0$. Then there exists $n$ exceptional line bundles $L_j$, $j=0,\ldots,n-1$, unique up to torsion by a flat line bundle  $F\in H^1(S,\bb C^\star)$ such that:
\begin{itemize}
 \item $E_{j}=c_{1}(L_{j})$, $0\leq j\leq n-1$ is a ${\mathbb Z}$-basis 
of $H^2(S,{\mathbb Z})$,
\item  $K_SL_j=-1$ and $L_{j}L_{k}=-\delta _{jk}$,
\item  $K_{S}=L_{0}+\cdots+L_{n-1}$ in 
$H^2(S,\bb Z)$
\item For every $i=0,\ldots,n-1$ and for every flat line bundle $F\in H^1(S,\bb C^\star)$, 
$$\cal X(L_i\otimes F)=0.$$
\item If $h^0(S,L_i\otimes F)\neq 0$, there exists an exceptional divisor $C_i$ and a (perhaps trivial) flat effective divisor $P_i$ such that $L_i\otimes F=[C_i+P_i]$.
 \end{itemize}
  \end{Th}
  Proof:   1) By Donaldson theorem \cite{DON},there exists a $\bb Z$-basis $(E_i)_{i}$, $0\le i\le n-1$, of $H^2(S,\bb Z)/Torsion$ such that $E_iE_j=-\d_{i,j}$. Moreover, since $p_g=h^2(S,\cal O_S)=0$, the exponential exact sequence $0\to \bb Z\to\cal O\to \cal O^\star\to 0$ yields line bundles $L_i$ such that $E_i=c_1(L_i)$. The line bundles $L_i$ are unique up to tensor product by flat line bundles. For a surface of class VII, the group of flat line bundles is $H^1(S,\bb C^\star)$. Let $c=\sum n_iE_i\in H^2(S,\bb Z)$. Then $c^2=-\sum n_i^2$, therefore $c^2=-1$ \iff $c=\pm E_i$. Replacing if necessary $L_i$ by $L_i^{-1}$ we may suppose that $KL_i\le 0$ for $i=0,\ldots,n-1$. By Riemann-Roch formula
  $$\cal X(L_i\otimes F)=\cal X(\cal O_S)+\frac{1}{2}(L_i^2-KL_i)=\frac{1}{2}(-1-KL_i)\in\bb Z,\leqno{(\ast)}$$
  therefore $KL_i\le -1$. Since $(L_i)$ is a $\bb Z$-base of $H^2(S,\bb Z)$, $K=\sum_i k_iL_i$ with $k_i=-KL_i\ge 1$. From $-n=K^2=-\sum_{i=0}^{n-1} k_i^2$ we deduce that $k_i=1$ for $i=0,\ldots,n-1$. From $(\ast)$ we obtain
  $$\cal X(L_i\otimes F)=0.$$\\
 2)  If $h^0(S,L_i\otimes F)>0$, then $L_i\otimes F=[C_i]$ where $C_i$ is an effective divisor. Let
  $$C_i=n_1G_1+\cdots+n_pG_p$$
be a decomposition into irreducible components.
 \begin{itemize}
 \item If $G$ is an elliptic curve or a rational curve with a double point then $KG=-G^2\ge 0$.
 \item If $G$ is  a nonsingular rational curve, $KG=-2-G^2\ge -1$ and $KG=-1$ \iff $G$ is an exceptional curve of the first kind.
 \end{itemize}
 Therefore the condition 
 $$-1=KC_i=\sum_i n_iKG_i$$
  implies that there is an exceptional curve of the first kind, say $G_p$. Now we prove the result by induction on $n=b_2(S)\ge 1$.\\
 If $n=1$, there is only one exceptional line bundle $L$ and if $h^0(S,L\otimes F)\neq 0$, $S$ is not minimal, hence a blow-up of a surface $S'$ with $b_2(S')=0$. Then $L\otimes F=[C+P]$ where $C$ is an exceptional curve of the first kind and $P=0$ (if $S'$ has no curve) or $P$ is flat perhaps not trivial (if $S'$ is a Hopf surface). Suppose that $n>1$:
\begin{itemize}
\item   If $ G_p\sim E_i$, then $-F=[n_1G_1+\cdots+n_{p-1}G_{p-1}+(n_p-1)G_p]$ is a flat line bundle and $L'_i=L_i\otimes F=[G_p]$.
\item If $G_p\sim E_j$, $j\neq i$,  there is a flat line bundle $F$ such that $L'_j=L_j\otimes F=[G_p]$.
\end{itemize}
 Therefore, replacing $L_j$ by $L'_j$, and changing if necessary the numbering we may suppose that $L_{n-1}=[C_{n-1}]$ with $C_{n-1}$ an exceptional curve of the first kind. Let $\P:S\to S'$ be the blowing-down of $C_{n-1}$. Since $L_i$, $i\neq n-1$ is trivial in a neighbourhood of $C_{n-1}$, $L'_i=\P_\star L_i$ is a line bundle such that $L_i=\P^\star \P_\star L_i$ and we check easily, by projection formula, that $(L'_i)_{0\le i\le n-2}$ is a family of exceptional line bundles. If $h^0(L_i)\neq 0$, then $h^0(L'_i)\neq 0$, whence by induction hypothesis, $L'_i=[C'_i+P'_i]$ with $C'_i$ an exceptional divisor and $P'_i$ an effective flat divisor. Therefore $L_i=\P^\star L'_i=[\P^\star C'_i+\P^\star P'_i]$.
\hfill $\Box$\\

\subsection{Families of exceptional or flat line bundles}
 If $S$ is minimal then, for any $i=0,\ldots,n-1$ and for any flat line bundle $F$, $h^0(S,L_i\otimes F)=0$. In all known examples $S$ has a deformation into a non minimal one, hence we consider now versal deformation of $S$.\\
Let $\P:\cal S\to B$ be the versal deformation of $S\simeq S_0$, where $B$ is the unit ball of $\bb C^N$, $N=\dim H^1(S,\T)$.  Standard arguments of spectral sequences yield
\begin{Prop}\label{ExcepFamily} For $i=0,\ldots,n-1$, there exist line bundles $\p_i:\cal L_i\to \cal S$ such that for every $u\in B$, the restriction $\p_{i,u}:L_{i,u}\to S_u$ is an exceptional line bundle,\\
Moreover, if $\cal K\to \cal S$ is the relative canonical bundle, we have
$$\cal K\sim \bigotimes_{i=0}^{n-1}\cal L_i$$
in $H^2(\cal S,\bb Z)\simeq H^2(S_0,\bb Z)$.
\end{Prop}

Since $\bb C^\star$ is commutative any representation $\rho:\p_1(S)\to GL(1,\bb C)\simeq \bb C^\star$ factorize through $H_1(S,\bb Z)$, therefore any representation (hence any flat line bundle), is defined by $\g\mapsto f$ with $f\in\bb C^\star$.  We shall denote by $F^f\in H^1(S,\bb C^\star)\simeq \bb C^\star$ this line bundle and we have defined a group morphism
$$\begin{array}{cccc}
\f:&\bb C^\star &\to & H^1(S,\bb C^\star)\\
&\l&\mapsto & F^\l
\end{array}$$

\begin{Lem}\label{tautbundle} 1) For any holomorphic function $f:B\to \bb C^\star$ there exists a unique flat line bundle $\cal F^f\to \cal S$ such that $(\cal F^f)_{\mid S_u}=F^{f(u)}$.\\
2) There exists over $\cal S\times \bb C^\star$ a flat line bundle, called the {\bf tautological flat line bundle} $\cal F$ such that for any $(u,\l)\in B\times\bb C^\star$, $\cal F_{S_u\times\{\l\}}=F^\l$. 
\end{Lem}
Proof: Let $\o:\tilde {\cal S}\to \cal S$ be the family of universal covering spaces of $\cal S$. Then the fundamental group $\p_1(S)$ operates diagonally on $\tilde{\cal S}\times\bb C$ by $\g.(p,z)=(\g.p,f(\P\o(p))z)$. The quotient manifold is $\cal F^f$.\\ 
A similar construction gives the tautological flat line bundle $\cal F$.\hfill $\Box$

\begin{Ex}\label{examplesgamma} {\rm
1) Suppose that  $S= S(F)$  is a primary Hopf surface defined by
$$F(z)=(\alpha_{1}z_{1}+sz_{2}^{m},\alpha_{2}z_{2}), \qquad 0<|\alpha_1|\leq
|\alpha_2|<1,\quad (\alpha_2^m-\alpha_1)s=0.$$
If $s=0$ (resp. $s\neq 0$), $S$ contains at least two elliptic curves $E_1$, $E_2$ (resp. only one elliptic curve $E_2$), where
$$E_{1}=\{z\neq 0\mid z_{1}=0\}/\{\alpha 
_{2}^{p}|p  \in{\bb Z}\}, \quad E_2=\{z\neq 0\mid z_{2}=0\}/\{\alpha _{1}^{p}|p\in{\bb Z}\}.$$
  Then for
$i=1,  2$, $\f
(\alpha _{i})=[E_{i}]$. In fact, if $S$ is a diagonal Hopf surface, the cocycle of the line 
bundle associated to 
$[E_{i}]$ is given by
$$(z_{1}, z_{2}, \lambda )\sim (\alpha _{1}z_{1}, \alpha _{2}z_{2}, 
\alpha _{i}\lambda ).$$\\ 
2)  Following \cite{D1}, let  $S= S(F)$ be the  minimal surface containing a GSS with 
$b_{2}(S)>0$ defined by
$$F(z)=(t^nz_1z_2^n+\sum_{i=0}^{n-1}a_it^{i+1}z_2^{i+1},tz_2).$$  
Then $S$ contains a cycle of
rational curves  $\Gamma=D_0+\cdots+D_{n-1} $ such that $D_i^2=-2$, $i=0,\ldots,n-1$ and $\Gamma ^{2}=0$. 
Let $t=tr DF(0)\neq 0$ be the trace of the surface, then  $\f 
(t)=[\Gamma ]$. In fact the equation of $\Gamma$ is $z_2=0$.\\
 If in the expression of $F$, there is at least one index $i$ such that $a_i\neq 0$, $S$ has no elliptic curve. If $M(S)$ is the intersection matrix  of the rational curves then $\det M(S)=0$, hence the curves do not generate $H^2(S,\bb Z)$; for every $m\ge 1$, every $F$ flat,  $H^0(S, -mK+F)=0$, therefore there is no NAC divisor.\\
 If $F(z)=(t^nz_1z_2^n,tz_2)$, i.e. $S$ is a Inoue (parabolic) surface, $S$ contains an elliptic curve $E$ and $-K=[E+\G]$.\\
3) Following \cite{I}, let $S=S_M=\mathbb H\times\mathbb C/G_M$ with $M\in SL(3,\mathbb Z)$ a unimodular matrix with eigenvalues $\alpha$, $\beta$, $\bar \beta$ such that $\alpha>1$, $\beta\neq\bar\beta$. Denote by $(a_1,a_2,a_3)$ a real eigenvector associated to $\alpha$ and $(b_1,b_2,b_3)$ an eigenvector associated  to $\beta$. It can be easily checked that $(a_1,b_1)$, $(a_2,b_2)$ and $(a_3,b_3)$ are linearly independent over $\mathbb R$. Let $G_M$ generated by
$$g_0:(w,z)\mapsto (\alpha w,\beta z),$$
$$g_i:(w,z)\mapsto (w+a_i,z+b_i) \quad {\rm for}\quad i=1,2,3.$$
If $G$ is generated by $g_i$, $i=1,2,3$, $\omega=dw\wedge dz$ is invariant under $G$ hence yields a 2-form on $\mathbb H\times\mathbb C/G$. Moreover, $g_0^\star \omega=\alpha\beta\omega$, hence yields a non-vanishing twisted 2-form over $S_M$ and $K=F^{1/\alpha\beta}$. A line bundle has no section for there is no curve.}
\end{Ex}

\subsection{Degeneration of blown-up Hopf surfaces}
All surfaces containing GSS are degeneration of blown-up primary Hopf surfaces as it can be easily checked using contracting germs of mappings. We show in this section that if a surface can be deformed into a non minimal one then over a Zariski open set in the base of the versal deformation, there are blown-up Hopf surfaces.\\
We need a lemma comparing the versal deformation of a surface $S$ with the versal deformation of a blowing-up $S'$ of $S$.
\begin{Lem}\label{DefEclatement} Let $S$ be a compact complex surface of the $VII$-class (not
necessarily minimal),  let
$\Pi :S'\rightarrow S$ be the blowing-up of $S$ at the point $z_{0}$ and
$$e_{z_{0}}:H^{0}(S, \Theta )\rightarrow T_{z_{0}}S$$
be the evaluation of global vector fields at $z_{0}$. Then, if $\cal V$ is a covering of $S$ such that $H^1(S,\T)=H^1(\cal V,\T)$,\\
1) $h^{0}(S', \Theta' ) = dim\ Ker\ e_{z_{0}}$;  \\
2)  There exists a covering ${\cal V}'=(V'_{i})_{i\geq 0}$ of $S'$ such that $H^1(S',\T')=H^1(\cal V',\T')$ with the following properties:
\begin{description}
\item{i)} $V'_{0}$ is the inverse image by $\Pi$ of a ball $V_{0}$ centered at $z_{0}$, 
\item{ii)} $V'_{0}$ meets only one open subset of the covering, say $V'_{1}$ along a
spherical shell,
\item{iii)} For all $i\geq 1$, the restriction of $\Pi$ on $V'_{i}$ is an isomorphism on its
image $V_{i}$,
\item{iv)} the canonical mapping $\P^\star:H^{1}({\cal V},\Theta)\rightarrow H^{1}({\cal
V}',\Theta')$ is injective,
\item{v)} A base of  $H^{1}(S', \Theta' )$ may be obtained from  a
base of $H^{1}({\cal V},\Theta)$  by adding  cocycles induced on $V'_{01}$ by (at most two)
non-vanishing vector fields
$Z^{i}$ on $V_{0}$ such that the vectors
$Z^{i}(z_{0})$ generate a suplementary subspace of $\ Im\
e_{z_{0}}$ in $T_{z_{0}}$.
\end{description}
 In particular 
$h^{1}(S', \Theta' ) = h^{1}(S, \Theta ) + codim\ Im\
e_{z_{0}}$.

\end{Lem} 
Proof: 1)
is clear.\\
 2) Let
${\cal U}=(U_{i})_{i\geq 1}$ be a locally finite covering of
$S$ such that
$H^{1}(S,
\Theta )=H^{1}({\cal U},\Theta )$. It may be supposed that $z_{0}\in U_{i}$ for
$i=1, \ldots p$. Denote by
$U'_{0}\subset\subset U_{0}$ balls centered at $z_{0}$ such that if $i>p$,  then
$U_{0}\cap U_{i}=\emptyset$. Now,  if $V_{0}=U_{0}$, $V_{1}=U_{1}\backslash
\overline{U'_{0}}$ and
$V_{i}=U_{i}\backslash
\overline{U_{0}}$,  for $i>1$, there are three coverings of $S$, ${\cal U}$,  ${\cal
U}_{0}=(U_{i})_{i\geq 0}$ and ${\cal V}=(V_{i})_{i\geq 0}$ related by the relation
$$ {\cal V}\prec {\cal U}_{0}\prec {\cal U}.$$
The canonical mappings 
$$H^{1}({\cal U},\Theta )\rightarrow H^{1}({\cal U}_{0},\Theta ) \rightarrow H^{1}({\cal
V},\Theta )$$
are isomorphism. We define a covering ${\cal V}'=(V'_{i})$ of $S'$ by
$V_{i}'=\Pi^{-1}(V_{i})$. The canonical morphism $\P^\star$ is clearly injective and the evident mapping
 $$s:H^{1}({\cal V}',\Theta' ) \rightarrow H^{1}({\cal V},\Theta )$$
is clearly surjective. 
Let $\xi\not\in Im\ e_{z_{0}}$ and $\theta$ a vector field on $U_{0}$ such that
$\theta(z_{0})=\xi$. Define $\eta=(\eta_{jk})\in Z^{1}({\cal V},\Theta)$ by
$\eta_{01}=\theta$ and
$\eta_{jk}=0$ if $\{j,k\}\neq\{0,1\}$. If $\eta'=\Pi^{\star}(\eta)\in Z^{1}({\cal
V}',\Theta')$, $\eta$ and $\eta'$ are cocycles such that $s([\eta'])=[\eta] =0$, but
$[\eta']\neq 0$. In fact if there exist  vector fields $X'_{0}$ on $V'_{0}$ and $X'$ on
$S'\backslash \Pi^{-1}(U'_{0})$ such that $\theta'=X'-X'_{0}$ on $V'_{01}$, we have
$\theta = \Pi_{\star}X'-\Pi_{\star}X'_{0}$ on $V_{01}$. But since a vector field extends
inside a ball, $\xi=\theta(z_{0})=\Pi _{\star}X'(z_{0})$ which is a contradiction.
Therefore,
$dim\ Ker\ s\geq codim\ Im\ e_{z_{0}}$ and it yields $h^{1}({\cal V}', \Theta' ) \geq
h^{1}({\cal V},
\Theta ) + codim\ Im\ e_{z_{0}}$. Now by Riemann-Roch-Hirzebruch-Atiyah-Singer theorem
we have, since $S$ and $S'$ are $VII$-class surfaces,
$$h^{1}(S,\Theta)=h^{0}(S,\Theta)+2b_{2}(S) \mbox{\quad and\quad}
h^{1}(S',\Theta')=h^{0}(S',\Theta')+2b_{2}(S') .$$ 
Using $b_{2}(S')=b_{2}(S)+1$, we obtain by a)
$$h^{1}(S', \Theta' ) = h^{1}(S, \Theta ) + codim\ Im\ e_{z_{0}}$$
therefore
$$
\begin{array}{lllll}
h^{1}({\cal V}, \Theta ) + codim\ Im\ e_{z_{0}}&=& h^{1}(S, \Theta ) +
codim\ Im\ e_{z_{0}}& =& h^{1}(S', \Theta' ) \geq h^{1}({\cal V}', \Theta' )\\
& \geq&
h^{1}({\cal V},\Theta ) + codim\ Im\ e_{z_{0}}&&
\end{array}
$$
which completes the proof. \hfill $\Box$

\begin{Th}\label{Th1} Let $S$ be a surface of class VII$_0^+$, $n=b_2(S)\ge 1$ and $\cal S\to B$ its versal deformation. Assume that there is a point $v\in B$, $\l_i\in\bb C$, $i=0,\ldots,n-1$ such that
$$H^0(S_v,L_{i,v}\otimes F^{\l_i})\neq 0.$$
Then there exists
\begin{itemize}
\item Holomorphic functions $c_i:B\to \bb C$,
\item Flat families of exceptional divisors $\cal C_i$ over $B\setminus H_i$, where $H_i=\{c_i=0\}$,
\end{itemize}
such that
\begin{itemize}
\item For every $u\in B\setminus H_i$, $[C_{i,u}]=L_{i,u}\otimes F^{c_i(u)}$ and $F^{c_i(v)}= F^{\l_i}$,
\item $S_u$ is minimal \iff $u\in M:=\cap_{i=0}^{n-1}H_i$,
\item $S_u$ is a blown-up Hopf surface
\iff $u\in B\setminus \cup_{i=0}^{n-1}H_i$.
 \end{itemize}
\end{Th}
Proof: 1) The surface $S_v$ contains an exceptional curve of the first kind $C$. Changing if necessary the numbering we may suppose that $C=C_{n-1}$ and  $L_{n-1,v}\otimes F^{\l_{n-1}}=[C_{n-1}]$. By stability theorem of Kodaira \cite{KO1}, there exists an open neighbourhood of $v\in B$ and a flat family $\cal C_{n-1}$ over $V$ of exceptional curves of the first kind. \\  2) Let $pr_{1}:{\cal S}\times {\bb
C}^{\star}\rightarrow {\cal S}$ be the first projection and define 
the 
line bundle
${\cal M}_{n-1}$ over ${\cal S}\times {\bb C}^{\star}$ by
$${\cal M}_{n-1}:=pr_{1}^{\star}{\cal L}_{n-1}\otimes_{\cal O_{\cal S\times \bb C^\star}}\cal F
$$ where $\cal F $ is the tautological flat line bundle of Lemma 
\thesection.\ref{tautbundle};  by  $p=\p_{n-1}\times Id$,  the
sheaf ${\cal M}_{n-1}$ is flat over $B\times{\bb C}^{\star}$. For 
every $(u,\alpha )\in B\times{\bb C}^{\star}$,  we have $({\cal 
M}_{n-1})_{|S_{u}\times\{\alpha \}}\simeq
({{\cal L}_{n-1}})_{|S_{u}}\otimes  F ^{\alpha }$.  By the 
semi-continuity theorem of
Grauert,  and because surfaces have no nonconstant meromorphic functions,
$$Z_{n-1}:=\Bigl\{(u, \alpha )\in B\times{\bb C}^{\star}\mid 
h^{0}\bigl(S_{u}\times\{\alpha \},
({\cal M}_{n-1})_{S_{u}\times\{\alpha \}}\bigr)= 1\Bigr\}$$
is an analytic subset of $B\times{\bb C}^{\star}$ and the dimension 
of the intersection
$Z_{n-1}\cap V\times{\bb C}^{\star}$ is  $N=\dim B$ by 1). \\
3) Let
$Z'_{n-1}$ be the irreducible component of $Z_{n-1}$ such that ${\cal
M_{n-1}}_{|p^{-1}(Z'_{n-1}\cap V\times{\bb C}^{\star})}=[{\cal 
C}_{n-1}]$. We have a flat family of curves
$$p:\cal C_{n-1}\to Z'_{n-1}$$
such that for $(u,\a)\in Z'_{n-1}$, $[C_{n-1,(u,\a)}]=L_{n-1,u}\otimes F^\a$.
For $p=0,1$, the functions
$$\begin{array}{ccc}
Z'_{n-1}&\to& \bb N\\
(u,\a)&\mapsto&h^p(C_{n-1,(u,\a)}, \cal O_{C_{n-1,(u,\a)}})
\end{array}$$
are constant (see \cite{BPV} p 96). For $u\in V$ and $\a$ such that the section of $L_{n-1}\otimes F^\a$ vanishes on an exceptional curve of the first kind,
$$h^0(C_{n-1,u}, \cal O_{C_{n-1,u}})=1,\quad h^1(C_{n-1,u}, \cal O_{C_{n-1,u}})=0$$
therefore $h^0(C_{n-1,u}, \cal O_{C_{n-1,u}})=1$ and $h^1(C_{n-1,u}, \cal O_{C_{n-1,u}})=0$ everywhere, whence over each $(u,\a)\in Z'_{n-1}$ the analytic set is connected and does not contain any elliptic curve or cycle of rational curves.\\
Now, we show that the intersection
$Z'_{i}\cap (V\times{\bb C}^{\star})$ contains only one irreducible 
component: In fact,  if
$z=(u, \alpha )$ and
$z'=(u, \alpha ')$ are two points in $Z'_{n-1}$ over $u\in V$, then
$$h^0(S_u,L_{n-1,u}\otimes F^\a)=h^0(S_u,L_{n-1,u}\otimes F^{\a'})=1,$$
whence $F^{\a/\a'}$ has a meromorphic section and by \cite{N1} (2.10), $F^{\a/\a'}=[D]$ with $D=mE+nF$, where $E,F$ are elliptic curves or cycles of rational curves such that $E^2=F^2=0$. It means that for $L_{n-1,u}\otimes F^\a=[C_{n-1,u}]$ and $L_{n-1,u}\otimes F^{\a'}=[C'_{n-1,u}]$ we have
$$[C_{n-1,u}]=L_{n-1,u}\otimes F^\a=L_{n-1,u}\otimes F^{\a'}\otimes F^{\a/\a'}=[C'_{n-1,u}+D]$$
and $C_{n-1,u}$ would not be connected, a contradiction.\\
 As consequence  $Z'_{i}$ cannot 
accumulate on 
$$B\times\{0\}\cap V\times{\bb
P}^{1}({\bb C}), \quad {\rm or}\quad B\times\{\infty\}\cap 
V\times{\bb P}^{1}({\bb 
C}).$$
By Remmert-Stein theorem the closure 
$$G_{n-1}=\overline{Z'_{n-1}}\subset B\times {\bb P}^{1}({\bb
C})$$
 is an irreducible analytic set of codimension one. The 
restriction
$$pr_{1}:G_{n-1}\rightarrow B$$
 is proper;  therefore,  $pr_{1}(G_{n-1})$ 
is an analytic subset
of $B$. Since it contains the open set $V$, $pr_{1}(G_{n-1})=B$. Now , 
$(G_{n-1}, pr_{1}, B)$
is a ramified covering which has only one sheet over $V$, hence 
$G_{n-1}$ is the graph of
a holomorphic  mapping $c_{n-1}:B\rightarrow {\bb 
P}^{1}({\bb
C})$. Define
$H_{n-1}:=c_{n-1}^{-1}(0)\cup c_{n-1}^{-1}(\infty)$.   So far,  we have
considered the family of curves
${\cal C}_{n-1}$ over
$Z'_{n-1}$. From now on, we shall consider it over
$B_{n-1}:=B\backslash H_{n-1}$. \\
4) If for $i\neq n-1$, $L_{i,v}\otimes F^{\l_i}=[C_{i,v}]$ with $C_{i,v}$ reducible, it means that $C_i=\p^\star(C'_i)$ where $\p$ is a finite sequence of  blowing-ups  and  $C'_i$ is an exceptional curve of the first kind.  By lemma (\ref{DefEclatement}), there is a deformation in which the blown-up points are moved outside $C'_i$. Therefore changing the point $v$ in $V'\subset V$ we may suppose that $C_{i,v}$ is an exceptional curve of the first kind and apply 3).\\
5) We show now that  for all $u\in B':=B\backslash\cup H_{i}$,  $S_{u}$ is 
a  blown-up Hopf
surface:  By Iitaka theorem we may blow down the
exceptional divisors ${{\cal C}_i}_{\mid B'}$ over $B'$. Let $p:\mathcal S\to\mathcal S'$ be the canonical mapping and $\Pi':{\cal S}'\rightarrow B'$ the induced family. 
Using classification of complex surfaces of class VII$_0$ with $b_2(S)=0$ (see \cite{LT}), $S'_{u}$ is a  Hopf surface or a Inoue surface \cite{I} of type $S_M$, $S_{N,p,q,r;t}^{(+)}$ or $S^{(-)}_{N,p,q,r}$. We have to exclude the last three types: By \cite{DOT2} we may suppose that $S$ has no non-trivial global vector field, therefore $N=h^1(S,\Theta)=2n$ and restricting if necessary $B$ we have $h^1(S_u,\Theta_u)=2n$ for all $u\in B$. We denote by $\mathcal K'$ the relative canonical bundle over ${\cal S}'$. For every $u\in B'$ there exists a unique $f(u)\in\mathbb C^\star$ such that $-K'_u=F^{f(u)}$; by Grauert theorem, $f:B'\to\mathbb C^\star$ is holomorphic. The relative canonical bundle $\cal K$ over $\mathcal S$ satisfies 
$$\mathcal K=p^\star\mathcal K'\otimes \bigotimes_{i=0}^{n-1}[\mathcal C_i]=F^{f^{-1}}\otimes \bigotimes_{i=0}^{n-1}[\mathcal C_i].$$
If there exists a point $u\in B'$ such that $S'_u$ is of type  $S_{N,p,q,r;t}^{(+)}$ or $S^{(-)}_{N,p,q,r}$, then  in a neighbourhood of $u$, all surfaces have the same type by \cite{I} and by the theorem of Bombieri \cite{I} p280, the function $f$ is real valued hence constant, in particular globally defined. Therefore $\bigotimes_{i=0}^{n-1}[\mathcal C_i] = F^f\otimes \mathcal K$ extends over $B$, which is impossible. Remains the case $S_M$ which is rigid: all the surfaces $S'_u$ are isomorphic hence by  (\ref{examplesgamma}) Example 3, 
$-K'_u=F^{\alpha\beta}$, therefore as before $\bigotimes_i[\mathcal C_i]$ should extend which is impossible.\\
6)  The functions $c_{i}$ are ${\mathbb C}$-valued:  In fact,  
suppose that 
$c_{i}^{-1}(\infty)\neq\emptyset$;  let  $A$ be an irreducible 
component of
$c_{i}^{-1}(\infty)$. Since $c$ has values in $\mathbb C$, there 
exists an index $j$ such that $c_{j}$ vanishes along $A$. Let  $\mu$ 
(resp. $\nu $) be the order of the pole (resp. zero) of $c_{i}$
(resp. $c_{j}$) at a point $a\in A$. Then on a disc $\D$ containing $a$,
$g:=c_{i}^{\nu }c_{j}^{\mu}\in{\cal O}^{\star}(\D)$ and the family of 
positive 
divisors
$$[\nu {\cal C}_{i}+\mu {\cal C}_{j}]={\cal L}_{i}^{\nu }\otimes 
{\cal 
L}_{j}^{\mu } \otimes F^{g}$$
extends on $\D$. But it implies that ${\cal C}_{i}$ and 
${\cal 
C}_{j}$   extend also,  which
yields a contradiction. \hfill $\Box$

\begin{Def} Let $S$ be a surface of class VII$_0^+$. We shall say that $S$ is a degeneration of blown-up primary Hopf surfaces if there is a deformation $\cal S\to \D$ over the unit disc of $S\simeq S_0$, such that $S_u$ is a blown-up primary Hopf surface for $u\neq 0$.
\end{Def}

If $\p_1(S)=\bb Z$,  a surface which can be deformed into a $n$ times blown-up surface,  is a degeneration of blown-up primary Hopf surfaces. The surface $S_u$ is defined  by a contraction
$$F_{u}(z)=(\alpha_{1}(u)z_{1}+s(u)z_{2}^{m},\alpha_{2}(u)z_{2}) $$
$${\rm with}\quad 0< | \alpha_{1}(u)|\leq |\alpha_{2}(u)|<1 \ 
{\rm and}\ s(u)(\alpha_{2}^{m}(u)-\alpha_{1}(u))=0$$
or
$$F_{u}(z)=(\alpha_{1}(u)z_{1},\alpha_{2}(u)z_{2}+s(u)z_{1}^{m})$$
$$ {\rm with}\quad 0<| \alpha_{2}(u)|\leq |\alpha_{1}(u)|<1 \ 
{\rm and}\ s(u)(\alpha_{1}^{m}(u)-\alpha_{2}(u))=0.$$
  In both cases there is at least 
one elliptic 
curve $E_{2}$ (resp. $E_{1}$) induced by $\{z_{2}=0\}$  
(resp. $\{z_{1}=0\}$) and another $E_{1}$ (resp. $E_{2}$) induced 
by $\{z_{1}=0\}$ (resp. $\{z_{2}=0\}$) if $s=0$. The trace
$$t(u)=tr(S'_{u})=tr(S_{u})=tr(DF_{u}(0))=\alpha_{1}(u)+\alpha_{2}(u)$$
and the determinant
$$d(u)= \det DF_u(0)=\alpha_{1}(u)\alpha_{2}(u)$$
are bounded holomorphis functions on $B'$ hence extend on $B$. 
They depend only on the conjugation class of $F_{u}$.  We call $tr(u)$ the {\bf trace of the surface $S_u$}.\par
 By 
\cite{KO} II p696, we have the 
following 
description of the canonical bundle $K'_{u}$ of $S'_{u}$: If $S'_{u}$ 
is 
a diagonal Hopf surface then 
$$K'_{u}=[-E_{1,u}-E_{2,u}]=F^{(\alpha_{1}(u)\alpha_{2}(u))^{-1}}.$$
 If $s(u)\neq 0$, i.e. $S$ is not diagonal, then since 
$\alpha_{2}^{m}(u)=\alpha_{1}(u)$,
 
$$K'_{u}=[-(m+1)E_{2,u}]=F^{\alpha_{2}(u)^{-(m+1)}}= 
F^{(\alpha_{1}(u)\alpha_{2}(u))^{-1}}.$$

By \cite{BPV} I.9.1 (vii), the canonical bundle $K_{u}$ of $S_{u}$ 
over 
$B'$ satisfies
$$K_{u}=p^{\star}K'_{u}\otimes\bigotimes_{0\leq i\leq n-1}[C_{i,u}].$$
If  ${\mathcal K}$ denotes the relative canonical bundle and $ F$ 
 denotes the tautological flat line bundle,  
$${\mathcal K}= F^{(\alpha_{1}\alpha_{2})^{-1}}\otimes\bigotimes_{i}[{\mathcal C}_{i}] =F^{f}\otimes \bigotimes_{i}{\mathcal L}_{i}$$
where $$f:=\frac{\prod_{i}c_{i}}{\alpha_{1}\alpha_{2}}.$$
 Since $\mathcal K$ and $\mathcal L_{i}$ are globally defined on $B$, 
$F ^{f}$ 
 is globally defined  on
$B$ and $f\in {\mathcal O}^{\star}(B)$. Twisting, for example, $\cal 
L_{1}$ by $F^{f^{-1}}$, we may suppose that $f=1$ and then 
$$c=\prod_{i=0}^{n-1}c_{i}=\alpha_{1}\alpha_{2}=d\in {\mathcal 
O}^{\star}(B)$$
 satisfies $\| c\|_\infty\leq 1$.\\ \\
For contracting germs associated to surfaces containing GSS we refer to \cite{D1}.
\begin{Prop} Let $S$ be a surface obtained by degeneration of blown-up minimal 
Hopf surfaces. If $tr(S)\neq 0$, then $S$ contains a GSS and if $F:(\bb C^2,0)\to (\bb C^2,0)$ is a contracting germ of mappings associated to $S$, $tr(S)=tr DF(0)$.
\end{Prop}
Proof: One of the two functions $\a_i$, $i=1,2$, say $\a_2$, admit a limit $\a_2(0)\neq 0$, and $\a_1(0)=0$. Therefore $F^{\a_2}$ is globally defined. Since $|\a_1(u)|<|\a_2(u)|$, the elliptic curve $E_{2,u}$ induced by $z_2=0$ exists for all $u\not\in H$, hence $[E_{2,u}]=F^{\a_2}$. Since $u\mapsto h^0(S_u,F^{\a_2(u)})$ is upper-semicontinuous we obtain in $S$ a flat cycle of rational curves $E_{2,0}$. By Enoki theorem, $S$ contains a GSS.\hfill $\Box$

\section{Surfaces with a cycle of rational curves}
\begin{Def} Let $S$ be a surface and $C$ be an analytic subset  of dimension one. We say that $C$ is a $r$-cycle of rational curves if 
\begin{itemize}
\item $C$ is an elliptic curve when $r=0$, 
\item $C$ is a rational curve with a double point when $r=1$, and
\item $C=D_0+\cdots+D_{r-1}$,  is a reduced effective divisor such that $D_i$ is a non-singular rational curve for $i=0,\ldots,r-1$ and $D_0D_1=\cdots=D_{r-2}D_{r-1}=D_{r-1}D_0=1$, $D_iD_j=0$ in all other cases, when $r\ge 2$.
\end{itemize}
We denote by $\sharp(C)=r\ge 0$ the number of rational curves of the cycle $C$.
\end{Def}
With notations of theorem (\ref{Zbase}), we have
\begin{Th}[\cite{N2} (1.7)] \label{N2(1.7)}
Let $S$ be a VII$_{0}^+$ surface with $n=b_{2}(S)$. Assume that $S$ contains exactly one
 cycle $C$ of rational curves such that $C^2<0$. Then 
\begin{itemize}
\item  $C \sim -(L_{r}+\cdots+L_{n-1})$  for some $1\leq r\leq n-1$, if $S$  is not 
an odd Inoue-Hirzebruch surface, or
\item  $C \sim -(L_{0}+\cdots+L_{n-1})+F_{2}$, with $F_{2}$ of order two, 
 if $S$ is an odd Inoue-Hirzebruch surface.
\end{itemize}
\end{Th}
\begin{Rem} If $S$ is an even Inoue-Hirzebruch surface (=hyperbolic I-H surface), the two cycles $C$ and $C'$ satisfy
$$C\sim -(L_r+\cdots+L_{n-1}), \quad C'\sim -(L_0+\cdots+L_{r-1}).$$
\end{Rem}
The following lemma plays a crucial role in the computation of self-intersection of the cycle $C$.
\begin{Lem}[\cite{N2} (2.4)] \label{N2(2.4)}
Let $S$ be a surface of class VII$_0^+$ and without divisor $D$ such that $D^2=0$. Let $L_{I}:=\sum_{i\in I}L_{i}$, $L=L_{I}+F$, $F\in H^1(S,\mathcal 
O_{S}^\star)$ for a nonempty subset $I\subset[0,n-1]$. Then we have:\\
i) If $I\neq [0,n-1]$, then $H^q(S,L)=0$ for any $q$.\\
ii) If $L\otimes \mathcal O_{C}=\mathcal O_{C}$, then $I=[0,r-1]$, 
and $F=\mathcal O_{S}$, $K_{S}-L+C=\mathcal O_{S}$.\\
iii) If $LC_{i}=0$ for any irreducible component $C_{i}$ of $C$, then 
$I=[0,r-1]$.
\end{Lem}
\subsection{Surfaces with numerically m-anticanonical divisor}
\begin{Def}\label{DefNAC} Let $S$ be a compact complex surface of the VII$_0$ class  
with $b_2(S)>0$ and let
 $m\geq 0$ be an integer. 
We shall say that $S$ admits a {\bf \boldmath numerically $ m$-anticanonical
divisor}  if there exists  a divisor $D_m$,  and a flat line 
bundle $F$ in $H^1(S,\mathbb C^\star)$ such that in $H^1(S,\cal O^\star)$,
$$mK+F+[D_m]=0.$$
 We shall say that $S$ admits a numerially anticanonical divisor, shortly a {\bf NAC divisor}, if there exists an integer $m$ such that there exists a numerically $m$-anticanonical divisor.
\end{Def}
\begin{Rem}{\rm Let $S$ be a compact complex surface of the VII$_0$ class  
with $b_2(S)>0$. Let
 $m\geq 1$ be an integer and $F\in H^1(S,\mathbb C^\star)$. Then $H^0(S,mK\otimes F)=0$.
In fact, a section of $mK\otimes F$ must vanish; let $[\Delta]=mK\otimes F$ be the associated divisor. Then
$$0=(mK+F-\Delta)\Delta=mK\Delta-\Delta^2.$$
Since $mK\Delta\geq 0$ and $\Delta^2\leq 0$, $\Delta^2= 0$, therefore $b_2(S)=-K^2=0$.\\
This means that there is no numerically m-canonical divisors.}
\end{Rem}
For the convenience of the reader we recall basic known facts with slightly different proofs (see \cite{N2} lemma (3.1))
\begin{Lem} \label{CompConn} 1) If  $D_m$ exists, it is a positive divisor, hence $H^0(S, -mK-F)\neq 0$.\\
2) If an irreducible curve $C$ meets the support $|D_m| $ of $D_m$ then $C$ is contained in 
$|D_m|$.
\end{Lem}
Proof: 1) We denote by $D_i$ the irreducible components of $D_m$. Let $D_m=\sum_i k_i D_i=A-B$ where 
$A=\sum_{i\mid k_i>0}k_iD_i\geq 0$ and 
$B=\sum_{i\mid k_i<0}(-k_i)D_i \geq 0$ have no common component. 
If $B\neq 0$, $B^2= \sum (-k_i)BD_i<0$, there exists an index $j$ such that $BD_j<0$. Therefore, 
$$0\leq mKD_j=(-F-D_m)D_j=-D_mD_j=-AD_j+BD_j<0$$
\ldots a contradiction.\\
2) If $C$ is an irreducible curve and meets $ | D_m|$, $D_mC=-mKC\leq 0$, therefore 
$C$ is contained in 
$ |D_m |$.
\hfill $\Box$\\

If there exists a non trivial  divisor $A$ such that $A^2=0$, then a 
numerically m-anticanonical divisor exists if and 
only if $S$ is a Inoue surface. In this case, $S$ contains an 
elliptic curve $E$, a cycle $\Gamma$ of 
rational curves, $m=1$ and $K+E+\Gamma=0$.\\
When there is no non-trivial flat 
divisor, the  numerically m-anticanonical
divisor $D_m$ is clearly unique. \\

\subsection{The reduction lemma}
\begin{Def} The least integer $m$ such that there exists a numerically 
m-anticanonical divisor, is called the {\bf index }
of the surface $S$ and will be denoted by $m(S)$.\\
The index $m(S)$ of the 
surface $S$ is the lcm of the denominators of the coefficients 
$k_i$ of $D_1=\sum_i k_i D_i$.
\end{Def}

The proof of \cite{DOT3} (1.3) works under the following relaxed hypothesis: 
\begin{Lem}[Reduction lemma] \label{RevRam} Let $S$ be a surface of class  VII$_0$ 
with $b_2(S)>0$ and index $m=m(S)>1$.  Then there exists a 
diagram
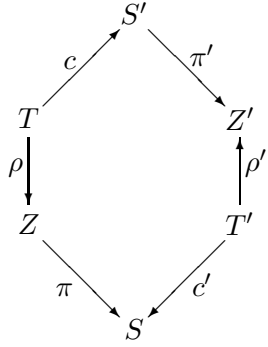
\DIAGV{50}
{}\n{}\n{S'}\n{}\n{}\nn
{}\n{\Near{c}}\n{}\n{\Sear{\pi'}}\n{}\nn
{T}\n{}\n{}\n{}\n{Z'}\nn
{\Sar{\rho}}\n{}\n{}\n{}\n{\naR{\rho'}}\nn
{Z}\n{}\n{}\n{}\n{T'}\nn
{}\n{\seaR{\pi}}\n{}\n{\swaR{c'}}\nn
{}\n{}\n{S}\n{}\n{}
\diag
where
\begin{description}
\item{i)} $(Z,\pi,S)$ is 
a m-fold  ramified covering space  of $S$, branched over $D_m$,  endowed with an automorphism 
group isomorphic to $\mathbb U_m$ which acts transitively on the 
fibers,
\item{ii)}  $(T,\rho,Z)$ is
 the minimal desingularization of $Z$,
\item{iii)} $(T,c,S')$ is the contraction of the (possible) exceptional 
curves of the first kind,
\item{iv)} $S'$ is a surface of  class VII$_0$,  with $b_2(S)>0$, with 
action of $\mathbb U_m$, with index $m(S')=1$,
\item{v)} $(S',\pi',Z')$ is the quotient space of $S'$ by $\mathbb 
U_m$, 
\item{vi)} $(T',\rho',Z')$ the minimal desingularization of $Z'$,
\item{vii)} $(T',c',S)$ is the contraction of the (possible) exceptional 
curves of the first kind,
\end{description}
such that the restriction over $S\setminus D$ is commutative, i.e.
$$\theta:=\pi\circ \rho\circ 
c^{-1}=c'\circ{\rho'}^{-1}\circ\pi':S'\setminus D'\to S\setminus D$$
and  $(S'\setminus D',S\setminus D)$ is a m-fold non ramified 
covering. 
Moreover 
\begin{itemize}
\item  $S$ contains a GSS  if and only if $S'$ contains a 
GSS, $S$ and
\item The maximal divisors $D$ and $D'$ of $S$ and $S'$ respectively have the same number of cycles and branches. 
\end{itemize}
\end{Lem}

\begin{Cor}\label{DmCycle} If $S$ admits a NAC divisor $D_m$, then the support of $D_m$ contains a cycle.
\end{Cor}
Proof: We may suppose that there is no divisor such that $D^2=0$. By reduction lemma it is sufficient to prove that the support of $D_{-K}$ contains a cycle.
Suppose that there exists a divisor $D_{-K}$ and a flat line bundle $F$ such that $K+D_{-K}+F=0$. 
By Cartan-Serre duality
$$h^2(S,O_S(-D_{-K}))=h^2(S,K+F)=h^0(S,-F)=\left\{
\begin{array}{cl}0&{\rm if\;F\neq 0}\\1&{\rm if\;F= 0}\end{array}
\right.$$
By Riemann-Roch formula,
$$h^0(S,-F)-h^1(S,-F)=h^0(S,-F)-h^1(S,-F)+h^0(S,K+F)=\chi(F)=0,$$
hence 
$$h^1(S,-F)=\left\{
\begin{array}{cl}0&{\rm if\;F\neq 0}\\1&{\rm if\;F= 0}\end{array}
\right.$$
We have
$$0\to \mathcal O_S(-D_{-K})\to \mathcal O_S \to \mathcal O_{D_{-K}}\to 0$$
\begin{itemize}
\item If $F\neq 0$, the long exact sequence yields
$$h^0(S,\mathcal O_{D_{-K}})=h^1(S,\mathcal O_{D_{-K}})=1.$$
The support of $D_{-K}$ is connected and by \cite{N1} (2.7), $h^1(S,\mathcal O_{(D_{-K})_{red}})\geq 1$
hence $(D_{-K})_{red}$ contains a cycle of rational curves.
\item If $F=0$, the associated long exact sequence and \cite{N1} (2.2.1) imply
$$1\leq h^0(S,\mathcal O_{D_{-K}})=h^1(S,\mathcal O_{D_{-K}})\leq 2.$$
As before $h^1(S,\mathcal O_{(D_{-K})_{red}})\geq 1$ and there is at least one cycle.\\
If $h^1(S,\mathcal O_{(D_{-K})_{red}})=2$, then
by the already quoted result (2.2.1) there are two cycles of rational curves.
\end{itemize}
\hfill $\Box$

\subsection{Characterization of Inoue-Hirzebruch surfaces}

\begin{Lem}\label{MinDK} Let $S$ be a surface with a NAC divisor $D_{-K}=\sum k_iD_i$. We suppose that the maximal divisor contains a cycle of rational curves $C=D_0+\cdots+D_{s-1}$ with $s\ge 1$ irreducible curves. If there exists $j\le s-1$ such that $k_j=1$, then $k_i=1$ for all $i=1,\ldots,s-1$ and $C$ has no branch.
\end{Lem}
Proof: {\bf\boldmath Case $s=1$}: Since $D_0$ is a rational curve with a double point, the adjunction formula yields
$$D_0^2=-KD_0=D_{-K}D_0=\sum k_iD_iD_0=D_0^2+\sum_{i>0}k_iD_iD_0.$$
hence $\sum_{i>0}k_iD_iD_0=0$ and $C$ has no branch.\\
 {\bf\boldmath Case $s\ge 2$}: By adjunction formula,
$$2+D_j^2=-KD_j=\sum_i k_iD_iD_j=D_j^2 + \sum_{i\neq j}k_iD_iD_j,$$
whence $2=\sum_{i\neq j}k_iD_iD_j$. Since $\sum_{i\neq j}D_iD_j\ge 2$, $D_j$ meets at most two curves $D_{j-1}$ and $D_{j+1}$ (one if the cycle contains two curves) and  $k_{j-1}=k_{j+1}=1$. By connectivity we conclude.\hfill $\Box$

\begin{Prop}\label{DivAntiCan} Let $S$ be a surface with a NAC divisor $D_{-K}=\sum k_iD_i$ and let $C=D_{0}+\cdots +D_{s-1}$ be a cycle contained in the support of $D_{-K}$. One of the following conditions holds:
\begin{description}
\item{i)} There exists an index $0\le j\le s-1$ such that $k_j=1$, then $S$ is a Inoue surface or a Inoue-Hirzebruch surface.
\item{ii)}  For every $0\le j\le s-1$, $k_j\ge 2$, then $C$ has at least one branch and the support 
$D=|D_{-K}|$ of $D_{-K}$ is connected. More precisely, if 
$k=\max \{k_{i}\mid 0\leq i\leq s-1\}$, there 
exists a curve $D_{j}$ and a branch $H_{j}>0$ 
such that $k_{j}=k$ and 
$D_{j}H_{j}>0$.
\end{description}
In particular, each connected component of $| D_{K}|$ contains 
a cycle.
\end{Prop}
\begin{Def} With preceeding notations, a curve $D_{j}$ in the cycle 
such that $D_{j}H>0$ will be called the root of the branch $H$.
\end{Def}
Proof of (\ref{DivAntiCan}):  Taking if necessary a double covering it may be 
supposed, by 
\cite{N1}(2.14), that 
the cycle  has at least two curves, hence all the curves are 
regular. Let 
$$D_{-K}=A+B=\sum_{i=0}^{p-1}k_{i}D_{}+\sum_{i=p}^{p+q}k_{i}D_{i}$$
where the support of $A=\sum_{i=0}^{p-1}k_{i}D_{i}$ is the connected 
component of the cycle $C=\sum_{i=0}^{s-1}D_{i}$.\\
If $B$ contains another cycle resp. an elliptic curve), then $S$ is a Inoue-Hirzebruch 
surface \cite{N1} (8.1) (resp. a Inoue surface \cite{N1} (10.2)); suppose therefore that the support $|B|$ of $B$ is 
simply-connected. We have to show that $|D_{-K}|$ is connected: In fact, let $B_{0}$ be a connected component of 
$|B|$. There is a proper mapping $p:S\to S'$ onto  a normal
 surface $S'$ with normal singularities $a=p(|A|)$ and $b=p(B_{0})$. 
Since $B$ 
 is simply connected,
  $F$ is trivial on a strictly 
 pseudo-convex neighbourhood $U$ of 
 $B$ and thus a holomorphic section of $-K-F$ yields a non vanishing 
 holomorphic 2-form on $U\setminus B$, i.e. $(S',b)$ is Gorenstein. 
If $(S',b)$ would be an elliptic singularity, a two-fold covering $T$ 
of $S$ should contain three  
 exceptional connected divisors such that their contractions $q:T\to 
 T'$ would fullfil $h^0(T',R^1q_{\star}\mathcal O_{T})=3$. 
However, by Leray spectral 
 sequence, there is an exact sequence
 $$0\to H^1(T',\mathcal O_{T'})\to H^1(T,\mathcal O_{T})\to 
 H^0(T',R^1q_{\star}\mathcal O_{T})\to H^2(T',\mathcal O_{T'})\to 
 H^2(T,\mathcal O_{S})\leqno{(\dag)}$$
 where 
 $$p_{g}=h^2(T,\mathcal O_{T})=0,\quad  q=h^1(T,\mathcal 
 O_{T})=1\quad {\rm and}\quad h^2(T',\mathcal O_{T'})=h^0(T',\omega 
 _{T'})\leq 1,$$
 by Serre-Grothendieck duality and because $T'$ has no non-constant meromorphic functions. 
 By ($\dag$) we obtain a contradiction, hence 
  $(S',b)$ is a Gorenstein rational 
 singularity, hence a Du Val singularity. However, such a singularity has a trivial canonical 
divisor, 
 therefore $B=0$ and $D_{-K}=A$.  \\

i) Suppose that there is an index $j\le s-1$ such that $k_j=1$. By lemma (\ref{MinDK}), $C$ has no branch and $s=p$. By adjunction formula, we have for every $0\leq i\leq s-1$,
$$D_{i}^2+2=-KD_{i}=D_{-K}D_{i}=k_{i-1}+k_{i+1}+k_{i}D_{i}^2,$$
therefore
$$(k_{i-1}-1)+(k_{i+1}-1)=(k_{i}-1)(-D_{i}^2). \leqno{(\star)}$$
 By $(\ast)$,  $k_{i}=1$ for all $0\leq i\leq p-1$.
\begin{itemize}
 \item If for every $0\le i\le s-1$, $D_i^2=-2$,  $S$ is a Enoki surface by \cite{E}. Moreover, by hypothesis, $S$ admits a NAC divisor, hence is a Inoue surface or
 \item There is at least one index $k\le s-1$ such that $D_{k}^2\leq -3$ whence $C^2<0$. If there is no other cycle, we have shown that $D_{-K}=C$, hence
$$-b_{2}(S)=K^2=C^2,$$
$S$ is  a Inoue-Hirzebruch surface by \cite{N1} (9.2).
 \end{itemize}

 In case ii), if the maximum $k$ is reached at a curve $D_{i}$, $i\le s-1$, which is not a root,  then applying again ($\star$),  we obtain 
 $k_{i-1}=k_{i+1}=k_{i}=k$. By connexity, we reach a root.\hfill $\square$

\begin{Prop}\label{arbre-IH} Let $S$ be a surface of class VII$_0^+$ admitting a NAC divisor $D_m$. Then one of the following conditions is fulfilled:\\
 i) The maximal divisor $D$ is connected and contains a cycle with at least one branch;\\
 ii)  $S$ is a Inoue surface or a (even or odd) Inoue-Hirzebruch surface.
 \end{Prop}
 Proof: By (\ref{DmCycle}) $S$ contains at least one cycle
 By (\ref{RevRam}) $S'$ admits a NAC divisor $D_{-K}$. Then by (\ref{DivAntiCan}), and (\ref{RevRam}) $S'$ and $S$ are of the same type.\hfill $\Box$

\subsection{Surfaces with singular rational curve}

\begin{Th}\label{Singrat} Let $S$ be a surface with $n=b_2(S)$ containing a singular 
rational curve $D_0$ with a double point. Then \\
1) $D_{0}\sim -(L_{1}+\cdots+L_{n-1})$ and the connected component 
containing $D_{0}$ is $D_{0}+D_{1}+\cdots+ D_{p}$ for $1\leq p\leq n-1$, 
$D_{i}\sim L_{i}-L_{i-1}$, $1\leq i\leq p$. In particular 
$D_{0}^2=-(n-1)$, $D_{1}^2=\cdots=D_{p}^2=-2$.\\
2) The following conditions are equivalent:\\
i) There exists an integer 
$m\geq 1$ such that there exists a numerically $m$-anticanonical 
divisor,\\
ii) $S$ contains a GSS.\\
When these conditions are fulfilled then $S$ is either a 
Inoue-Hirzebruch surface and $m(S)=1$ or 
$D_0$ has a branch and its index satisfies $m(S)=n-1$.
\end{Th}
Proof: 1) By (\ref{N2(1.7)}),  there 
exists $r\geq 1$ such that $D_{0}=-(L_{r}+\cdots+L_{n-1})$. We have 
$L_{0}D_{0}=0$, hence $r=1$ by (\ref{N2(2.4)}) iii). We show by induction on 
$i\geq 1$ that if there is a non singular rational curve $D_{i}$ such that 
$(D_{0}+\cdots+D_{i-1})D_{i}\neq 0$ then $D_{i}$ is unique and 
$D_{i}=L_{i}-L_{i-1}$. By \cite{N1} (2.2.4), and unicity
$D_{0}D_{i}=\cdots =D_{i-2}D_{i}=0$, $D_{i-1}D_{i}=1$ if $i\ge 2$.\\
Suppose that there is a curve, say $D_1$ such that $D_{0}D_{1}=1$. 
By \cite{N2} (2.6) (see lemma \ref{N2(2.5)} below) we set $D_{1}\sim L_1-L_{I_1}$, 
then 
$$1=D_{0}D_{1}=-(L_1+\cdots+L_{n-1})(L_1-L_{I_1})=1+(L_1+\cdots+L_{n-1})L_{I_1},$$
hence $I_1=\{0\}$ and $D_1\sim L_1-L_0$.
 If there were another curve, say $D_2$ meeting $D_0$, 
the same argument shows that $D_2=L_2-L_0$, but in this case $D_1D_2=-1$ and it is impossible. It shows 
the unicity. Suppose 
that for $i> 1$, $D_i$ exists and $D_{i}\sim L_i-L_{I_i}$. 
We have the equations
$$\begin{array}{cclclcl}
0&=&D_{0}D_{i}&=&-(L_{1}+\cdots+L_{n-1})(L_i-L_{I_i}) &=&1+(L_{1}+\cdots+L_{n-1})L_{I_i}\\
0&=&D_{1}D_{i}&=&-(L_{1}-L_{0})(L_i-L_{I_i}) &=&-L_1L_{I_i}+L_0L_{I_i}\\
&\vdots&&&&\vdots\\
0&=&D_{i-2}D_{i}&=&-(L_{i-2}-L_{i-3})(L_i-L_{I_i}) &=&-L_{i-2}L_{I_i}+L_{i-3}L_{I_i}\\
1&=&D_{i-1}D_{i}&=&-(L_{i-1}-L_{i-2})(L_i-L_{I_i}) &=&-L_{i-1}L_{I_i}+L_{i-2}L_{I_i}
\end{array}$$
which yield from the last to the second  equation
$$i-1\in I_i,\  i-2\not\in I_i,\,\ldots,\,0\not\in I_i.$$
The first one implies that $I_i$ contains exactly one index, hence $$D_{i}\sim L_{i}-L_{i-1}.$$
 We prove unicity as for $i=1$.\\

2) By 1) the intersection matrix of the connected component of the 
cycle is the $(p+1,p+1)$ matrix, $p\ge 0$,
$$M=\left(
\begin{array}{crrrrr}
-(n-1)&1&0&\cdots&\cdots&0\\
1    &-2&1&0&&\vdots\\
0&1&-2&\ddots&\ddots&\vdots\\
\vdots&\ddots&\ddots&\ddots&1&0\\
\vdots&&\ddots&1&-2&1\\
0&\cdots&\cdots&0&1&-2
\end{array}
\right)$$
and it is easy to check that
$$\det M=(-1)^{p+1}\bigl[(n-1)(p+1)-p\bigr].$$
A numerically m-anticanonical divisor $D_{m}=\sum_{i=0}^pk_iD_{i}$ 
supported by the connected component containing the cycle satisfies 
the linear system
$$M\left(\begin{array}{c}k_{0}\\ k_1\\ \vdots\\k_{p}\end{array}\right) = 
\left(\begin{array}{c}-m(n-1)\\ 0\\ \vdots\\0\end{array}\right).$$
Therefore
$$k_{i}=\frac{m(n-1)(p+1-i)}{(n-1)(p+1)-p},\quad 
{\rm for\;every\; } 0\leq i\leq p.\leqno{(\ast)}$$

Notice that 
$$k_{i}=(p+1-i)k_{p}\quad 
{\rm for\;every\; } 0\leq i\leq p,$$
therefore all $k_{i}$ are integers if and only if $k_{p}$ is an 
integer.
By assumption 
 $$D_{m}^2=(-mK)^2=-m^2n,$$
  and we have $D_mD_i=-mKD_i=0$ for $i=1,\ldots,p$.  We have to consider two cases:\\
  {\bf First case}: $p\ge 1$ then
$$-m^2n=\sum_{i=0}^p k_{i}D_{i}D_{m}=k_{0}D_{0}D_{m}=k_{0}^2 
D_{0}^2+k_{0}k_{1} = -\frac{m^2 (n-1)^2(p+1)}{(n-1)(p+1)-p}$$
and this condition is equivalent to
$$p=n-1.$$
Hence $S$ has a numerically m-anticanonical divisor \iff $S$ contains $n$ rational curves \iff $S$ contains a GSS by \cite{DOT3}.\\
Replacing in ($\ast$) we obtain
$$k_{i}=\frac{m(n-i)}{n-1},\quad 
{\rm for\;every\; } 0\leq i\leq n-1.$$
Hence we conclude 
$$m(S)=n-1.$$
{\bf Second case}: $p=0$, then
$$-m^2n=k_{0}^2 D_{0}^2=-k_0^2(n-1)$$
However, this equation has no solution in integers, therefore the support of $D_m$ is not connected. We conclude by (\ref{arbre-IH}).

\hfill $\square$

\begin{Cor} If $S$ in the VII$_0^+$ class has a NAC divisor and contains a rational curve with a double point, then $S$ contains a GSS.
\end{Cor}

\subsection{Cycles with at least two rational curves}
\begin{Lem} [\cite{N2} (2.5), (2.6), (2.7)] \label{N2(2.5)}
 1) Let $D$ be a nonsingular rational curve.  If  $D\sim a_0L_0+\cdots+a_{n-1}L_{n-1}$, then there exists a unique $i\in [0,n-1]$ such that $a_i=1\;{\rm or }\; -2$, and $a_j=0\;{\rm or}\; -1$ for $j\neq i$.\\
 2) If $D$ is not contained in the cycle then $D\sim L_i-L_I$ for some $i\in[0,n-1]$ and $I\subset[0,n-1]$, $i\not\in I$.\\
 3) Let $D_1$ and $D_2$ two distinct nonsingular rational curves such that $D_1\sim L_{i_1}-L_{I_1}$, $D_2\sim L_{i_2}-L_{I_2}$, then $i_1\neq i_2$.
 \end{Lem}

\begin{Prop}\label{FormeNormaleCourbe} Let $S$ be a surface containing a cycle $C=D_0+\cdots+D_{s-1}$ of $s\ge 2$ nonsingular rational curves. Suppose that there exists a rational curve $E$ such that $E.C=1$, and let $D=\sum_i D_i$ be the maximal connected divisor containing $C$.  Then
every curve $D_i$ is of the type $D_i\sim L_i-L_I$ for $i\not\in I$.
\end{Prop}
Proof:    By (\ref{N2(2.5)}), every curve $D_i$ is of type $D_i\sim L_i-L_I$ (type {\bf a}) or $D_i\sim -2L_i-L_I$ (type {\bf b}). Suppose that $C$ contains at least one curve of type {\bf b}. Since $(-2L_i-L_I)(-2L_j-L_J)\le 0$ two such curves cannot meet, in particular if $s=2$, $C$ contains at most one curve of type {\bf b}.\\
Suppose $s\ge 3$. If $D_i$ of type {\bf b}, $D_i$ meets curves $D_j$ of type {\bf a}. We have 
$$1=(-2L_i-L_I)(L_j-L_J)=-2L_iL_j- L_jL_I + 2L_iL_J + L_IL_J$$
and 
\begin{itemize}
\item either $j\in I$, $i\neq j$, $i\not\in J$, $I\cap J=\emptyset$, whence
$$D_i+D_j\sim -2L_i-L_I+L_j-L_J=-2L_i - L_{I\cup J\setminus\{j\}}.$$
\item or $j\not\in I$, $i=j$, $i\not\in J$ and $I\cap J$ contains one element, say $k$. Setting $I'=I\setminus\{k\}$ and $J'=J\setminus\{k\}$ we obtain
$$D_i+D_j\sim -2L_k - L_{\{i\}\cup I'\cup J'}.$$
\end{itemize}
 By (\cite{N2} (1.7)), there is a proper smooth family of compact surfaces $\p:\cal S\to \D$ over the unit disc, a flat divisor $\cal C$ such that $C_0=C$ and for $u\neq 0$, $C_u\sim D'+(C-D_i-D_j)$, with $D'\sim D_i+D_j$ of type {\bf b}, in particular $\sharp(C_u)=\sharp(C)-1$. Therefore, repeting if necessary such a deformation we obtain a contradiction if $C$ contains two curves of type {\bf b}. Now if $C$ contains a curve of type {\bf b}, we take a double covering $p:S'\to S$ of $S$. The surface $S'$ of type VII$_0$ satisfies $b_2(S')=2b_2(S)$ and by \cite{N1} (2.14) contains a cycle of $2s$ rational curves. Applying the same arguments to $S'$, $p^\star D_i\sim -2p^\star L_i-p^\star L_I$ is the union of two rational curves of type {\bf b}, hence a contradiction. Finally, there is no curve of type {\bf b}.\hfill $\Box$

\begin{Lem} \label{DjDkDl} Let $S$ be a surface containing a cycle $C=D_0+\cdots+D_{s-1}$ of $s\ge 2$ nonsingular rational curves and $D=\sum_{i=0}^p D_i$ the maximal connected divisor containing $C$.\\
1) Let $D_j$, $D_k$ and $D_l$ three distinct nonsingular rational curves, 
$D_j\sim L_{i_j}-L_{I_j}$, $D_k\sim L_{i_k}-L_{I_k}$, $D_l\sim L_{i_l}-L_{I_l}$. Then $I_j\cap I_k \cap I_l=\emptyset$.\\
2) If $D_j$, $D_k$ are two distinct nonsingular rational curves, then $\Card (I_j\cap I_k)\le 1$.
\end{Lem}
Proof: The surface $S$ can be deformed into a blown-up  Hopf surface $\p:\cal S\to \D$ with a flat family $\cal C$, where $C_0=C=D_0+\ldots+D_{s-1}\sim -(L_r+\cdots+ L_{n-1})$, and $C_u$, $u\neq 0$,  is an  elliptic curve blown-up $n-r$ times. We denote by $\P=\P_{0,u}\cdots \P_{n-1,u}:S_u\to S'_u$ the composition of  blowing-ups and $(\P_{i+1}\cdots\P_{n-1})^\star(E_{i,u})\sim L_i$.  We have $D_j\sim L_{i_j}-L_{I_j}$, $j=0,\ldots,s-1$, hence $D_j$ is homologous to an exceptional rational curve of the first kind $E_{i_j}$ blown-up $\Card(I_j)$ times. Since an exceptional curve cannot blow-up three rational curves we have the first assertion. Moreover two distinct exceptional rational curve of the kind cannot be blown-up two times by the same curves, whence the second assertion. \hfill $\Box$

\begin{Prop}\label{sharpCb2} Let $S$ be a minimal surface containing a cycle $C=D_0+\cdots+D_{s-1}\sim -(L_r+\cdots+L_{n-1})$ of $s\ge 1$ nonsingular rational curves. Then, $r=s$ and numbering properly the line bundles $L_i$, for $i=0,\ldots,r-1$, we have $D_i\sim L_i-L_{I_i}$. Moreover $\cup_{i=0}^{s-1}I_i=[0,n-1]$ and
$$\sharp(C)-C^2=b_2(S).$$
\end{Prop}
Proof:  1) The case $s=1$ has been proved in \cite{N2} (see lemma (\ref{Singrat})).\\
2) {\bf \boldmath If $s=2$},  $D_0\sim L_{i_0}-L_{I_0}$, $D_1\sim L_{i_1}-L_{I_1}$. We have
$$2=D_0D_1=-L_{i_0}L_{I_1} -L_{i_1}L_{I_0}+ L_{I_0}L_{I_1},$$
whence $i_1\in I_0$, $i_0\in I_1$, $I_0\cap I_1=\emptyset$. Setting $I'_0=I_0\setminus\{i_1\}$ and $I'_1=I_1\setminus\{i_0\}$, we obtain
$$D_0\sim L_{i_0}-L_{i_1}-L_{I'_0}, \quad D_1\sim L_{i_1}-L_{i_0}-L_{I'_1}, \quad {\rm with}\; I'_0\cap I'_1=\emptyset, \; \{i_0,i_1\}\cap (I'_0\cup I'_1)=\emptyset.$$
Therefore
$$-(L_r+\cdots+L_{n-1})\sim C=D_0+D_1\sim -(L_{I'_0}+L_{I'_1})$$
i.e. $I'_0\cup I'_1=[r,n-1]$. Let $I=\{i_0,i_1\}\cup I'_0\cup I'_1$ and $I'=[0,n-1]\setminus I$. Of course, $L_{I'}.D_0=L_{I'}.D_1=0$, whence by (\ref{N2(2.4)}) 3), $I'=[0,r-1]$ which is impossible. Therefore $I'=\emptyset$ and $\{i_0,i_1\}\cup I'_0\cup I'_1=[0,n-1]$, i.e. $r=2$ and 
$$\sharp(C)-C^2=2+(n-2)=b_2(S).$$
3) {\bf \boldmath If $s\ge 3$}, $D_j\sim L_{i_j}-L_{I_j}$, $j=0,\ldots,s-1$ and we may suppose that $$D_0D_1=\ldots=D_{s-2}D_{s-1}=D_{s-1}D_0=1.$$
 Since there is a deformation of $S$ in which $D_j+D_{j+1}\sim L_{i_j}+L_{i_{j+1}} - L_{I_j} - L_{I_{j+1}}$ is deformed into a nonsingular rational curve $D'_u$ in a (perhaps non minimal) surface, contained in a cycle $C'_u$, hence $D_j+D_{j+1}$ must be of the homological type {\bf a},  and either $i_j\in I_{j+1}$, $i_{j+1}\not\in I_j$ or $i_j\not\in I_{j+1}$, $i_{j+1}\in I_j$. Moreover the equality
$$1=D_jD_{j+1}=-L_{i_j}L_{I_{j+1}} - L_{i_{j+1}}L_{I_{j}} + L_{I_j}L_{I_{j+1}}$$
implies that $I_j\cap I_{j+1}=\emptyset$, whence
$$\forall j, \quad 1=D_jD_{j+1}=-L_{i_j}L_{I_{j+1}} - L_{i_{j+1}}L_{I_{j}}.$$
4) Now, {\bf \boldmath we show by induction on $s\ge 2$ that 
for all $j$, $0\le j\le s-1$, there is a unique index, denoted $\s(j)$, such that $i_j\in I_{\s(j)}.$} This assertion is evident if $s=2$, therefore we suppose that $s\ge 3$. By 3) it is possible to choose the numbering in such a way that $i_{s-1}\in I_{s-2}$. We choose a deformation $\P:\cal S\to \D$ of $S$ over the disc endowed with a flat family $\cal C$ of curves such that $C_0=C$, and for $u\neq 0$,
$$D'_u\sim D_{s-2}+D_{s-1}\sim L_{i_{s-2}}+L_{i_{s-1}} - (L_{I_{s-2}}+ L_{I_{s-1}})=L_{i_{s-2}}- (L_{I'_{s-2}}+ L_{I_{s-1}}),$$
with $I'_{s-2}=I_{s-2}\setminus\{i_{s-1}\}$. Consider the cycle  $C_u=D_{0,u}+\cdots+D_{s-3,u} + D'_u$ of $s-1$ rational curves in the (perhaps non minimal) surface $S_u$; by the induction hypothesis, since $i_{s-1}\in I_{s-2}$ and $I'_{s-2}\cap I_{s-1}=\emptyset$, for all $0\le j\le s-2$, there is a unique index $\s(j)$, $0\le \s(j)\le s-1$ such that $i_j\in I_\s(j)$. Repeting this argument with $D_{s-3}+D_{s-2}$ we obtain the result.  Therefore we have a well-defined mapping $\s:\{0,\ldots,s-1\}\to \{0,\ldots,s-1\}$ such that for all $j$, $0\le j\le s-1$, $i_j\in I_{\s(j)}$.\\
5) Setting $I'_j=I_j\setminus \s^{-1}(j)$ we have
$$\begin{array}{lcl}
-(L_r+\cdots+L_{n-1})&\sim& C=D_0+\cdots +D_{s-1}\sim L_{i_0}+\cdots+L_{i_{s-1}} - (L_{I_0}+\cdots+L_{I_{s-1}})\\ 
&&\\
&\sim& - (L_{I'_0}+\cdots+L_{I'_{s-1}})
\end{array}$$
therefore
$$[r,n-1]=I'_0 \cup \cdots \cup I'_{s-1} \quad {\rm with}\quad \{i_0,\ldots,i_{s-1}\}\cap\Bigl( I'_0 \cup \cdots \cup I'_{s-1}\Bigr)=\emptyset$$
If $I=\{i_0,\ldots,i_{s-1}\}\cup I'_0 \cup \cdots \cup I'_{s-1}$ and $I'=[0,n-1]\setminus I$, we have $L_{I'}D_j=0$ for all $j=0,\ldots,s-1$, then if $I'\neq\emptyset$, lemma (\ref{N2(2.4)}) would imply that $I'=[0,r-1]$ and this is impossible, hence $I'=\emptyset$ and $\{i_0,\ldots,i_{s-1}\}\cup I'_0 \cup \cdots \cup I'_{s-1}=[0,n-1]$, i.e. $r=s$ and $\{i_0,\ldots,i_{r-1}\}=[0,r-1]$. Finally
$$\sharp(C)-C^2=s+(n-r)=n=b_2(S).$$
\hfill $\Box$

\begin{Cor} Let $S$ be a minimal surface with $b_2(S)\ge 1$. If $E$ is a $k$-cycle, $k=\sharp(E)\ge 0$, such that $i_\star H_1(E,\bb Z)=H_1(S,\bb Z)$, then $\sharp(E)-E^2=b_2(S)$.
\end{Cor}

\begin{Rem}{\rm If $[H_1(S,\bb Z):i_\star H_1(E,\bb Z)]=2$, $S$ is an odd Inoue-Hirzebruch surface with only one cycle $C$. Then we have $\sharp(C)-C^2=2b_2(S)$ (see \cite{N1} (2.13)).}
\end{Rem}

For surfaces containing a divisor $D$ such that $D^2=0$, i.e. Enoki surfaces, the situation is well understood: They all contain GSS, $D$ is a cycle of $b_2(S)$ rational curves $D=D_0+\cdots+D_{n-1}$ such that $D_0^2=\cdots=D_{n-1}^2=-2$ and $S$ admits a numerically anticanonical divisor \iff $S$ is a Inoue surface, in which case $h^0(S,-K)=1$.
For surfaces such that $D^2<0$ for any divisor we have the following theorem:

\begin{Th}\label{DNAC} Let $S$ be a compact complex surface of class VII$_0^+$. We suppose that there is no divisor such that $D^2=0$. Then the following properties are equivalent:
\begin{description}
\item{i)} $S$ contains a GSS,
\item{ii)} $S$ admits a NAC divisor.
\item{iii)} $S$ contains $b_2(S)$ rational curves,

\end{description}
\end{Th}
Proof:  $i) \Leftrightarrow iii)$ by \cite{DOT3}.\\
$i) \Rightarrow ii)$ The intersection matrix $M(S)$ is negative definite and the rational curves give a $\bb Q$-base of $H^2(S,\bb Q)$, whence there is an integer $m\ge 1$ and a divisor $D_m$ such that in $H^2(S,\bb Z)$, $mK+D_m=0$. By lemma (\ref{CompConn}), $D_m$ is effective. \\
$ii) \Rightarrow i)$ By lemma (\ref{RevRam}) we may suppose that the index satisfies $m(S)=1$, i.e. $S$ has a NAC divisor $D_{-K}$. By theorem (\ref{Singrat}) we may suppose that there is no singular curve. By lemma (\ref{DivAntiCan}) there are two cases: $S$ is a Inoue-Hirzebruch surface, in particular contains a GSS, or $S$ contains a cycle with at least one branch. Therefore we have to prove the result in the second case. Let $D=\sum_{i=0}^p D_i$ be the maximal divisor, where $C=D_0+\cdots+D_{s-1}$ is the cycle. By (\ref{FormeNormaleCourbe}), we have in $H^2(S,\bb Z)$,
$$-(L_0+\cdots+L_{n-1})=-K=\sum _{i=0}^p k_iD_i = \sum_{i=0}^p k_i(L_i-L_{I_i}),$$
where $k_i\ge 1$. If $p<n-1$, the curve $D_{n-1}$ is missing, hence there is exactly one index $j$ such that $n-1\in I_j$ and $k_j=1$. Moreover by Proposition (\ref{sharpCb2}), $\cup_{i=0}^{s-1}I_i=[0,n-1]$, therefore $j\le s-1$. By lemma (\ref{MinDK}), $S$ would be a Inoue-Hirzebruch surface, which is a contradiction.
\hfill $\Box$

\begin{Cor}[\cite{OTZ}] Let $S$ be a surface of class VII$_0$ with $b_2(S)>0$. Then $S$ is a Inoue-Hirzebruch surface if and only if  there exists  two twisted vector fields $\t_1\in H^0(S,\T\otimes F_1)$, $\t_2\in H^0(S,\T\otimes F_2)$, where $F_1$, $F_2$ are flat line bundles,   such that  $\t_1\wedge\t_2(p)\neq 0$ at at least one point $p\in S$.
\end{Cor}
Proof: $\t_1\wedge \t_2$ is a non trivial section of $-K\otimes F$ whence $S$ contains a GSS. We conclude by \cite{DO} th. 5.5.\hfill $\Box$

\section{On classification of bihermitian surfaces}
\subsection{Conformal and complex structures}
We consider  connected oriented conformal 4-manifolds $(M^4,c)$  with two complex (i.e. integrable almost-complex) structures $J_1$, $J_2$ which  induce the same orientation. Given a riemannian metric $g$ in $c$, $(M^4,g,J_1,J_2)$ is called a {\bf bihermitian surface} relatively to the conformal class $c$ if 
\begin{itemize} 
\item $J_i$ are orthogonal with respect to the metric, i.e. $g(J_iX,J_iY)=g(X,Y)$, $i=1,2$,
\item $J_1$ and $J_2$ are independent, i.e. there is a point $x\in M$ such that $J_1(x)\neq \pm J_2(x)$.
\end{itemize} 
The triple $(c,J_1,J_2)$ is called a bihermitian structure on $M^4$. 
If moreover $J_1(x)\neq \pm J_2(x)$ everywhere, the bihermitian structure $(c,J_1,J_2)$ is called {\bf stron\-gly bihermitian}. \\
Given two such almost-complex structures $J_1$ and $J_2$, denote by $f$ the smooth function, called the angle function
$$f=\frac{1}{4}\tr(^tJ_1J_2)=-\frac{1}{4} \tr(J_1J_2).$$
By Cauchy-Schwarz inequality, $|f|\leq 1$ and $f(x)=\pm 1$ \iff $J_1(x)=\pm J_2(x)$. Since $J_1J_2\in SO(4)$ and $J_2J_1$ is the inverse of $J_1J_2$, it is easy to check that
$$J_1J_2+J_2J_1=-2f\,Id.$$
Moreover $J_1(x)$ and $J_2(x)$ anticommute \iff $f(x)=0$.\\

Another conformal structure is provided by the Weyl curvature tensor: The riemannien curvature tensor $R$ of type (3,1) has a classical decomposition, under the orthogonal group $O(4)$, into three parts given by the scalar curvature, the Ricci curvature tensor without trace and the Weyl curvature tensor $W$ of type (3,1) which is a conformal invariant \cite{B}.\\
Let $\ast:\bigwedge^2T^\star M\to \bigwedge^2T^\star M$ be the star-Hodge operator, with $\star^2=Id$ and two eigenvalues $\pm 1$. We denote by $\bigwedge^2_+$ (resp. $\bigwedge^2_-$) the eigenspace associated to $+1$ (resp. $-1$). Since $\ast$ depends only on $c$, the splitting $\bigwedge^2T^\star M= \bigwedge^2_+\oplus \bigwedge^2_-$ is a conformal invariant. Let $W(g):\bigwedge^2T^\star M\to \bigwedge^2T^\star M$ be the Weyl curvature tensor of type (2,2), with restrictions $W_\pm\in End(\bigwedge^2_\pm)$ over $\bigwedge^2_\pm$. The riemannian conformal class of $(M^4,g)$  is called {\bf anti-self-dual (ASD)} if $\ast W(g)=-W(g)$, or equivalently $W_+=0$.\\
The geometric meaning of ASD condition stems from Atiyah-Hitchin-Singer theorem \cite{AHS}: Let 
$$Z=\{J\in SO(TM)\mid J^2=-Id\}=SO(TM)/U(2)\to M$$
be  the twistor space, i.e. the space of all orthogonal almost-complex structures over $(M,g)$, inducing the orientation of $M$. The fiber is isomorphic to the Riemann sphere $S^2$. Since there is only one complex structure on $S^2$, $Z$ is a differentiable fiber bundle with complex fiber $\bb P^1(\bb C)$. The complex structure on the fiber and on the base yield a canonical almost-complex structure $\bb J$ of the twistor space $Z$, which is not integrable in general. The theorem of Atiyah-Hitchin-Singer asserts that $\bb J$ is integrable \iff the metric $g$ is ASD, hence any compatible complex structure at a point $x\in M$ extends into a compatible complex structure over a neighbourhood of $x$. However, it does not extends to the whole manifold and perhaps $(M,c)$ admits no complex structure (for instance $S^4$). {\bf The aim is to classify compact 4-manifolds with several compatible complex structures, or at least to give necessary conditions for their existence}. When there are more than two compatible almost structures, we have
\begin{Prop}[\cite{PON}] If an oriented riemannian 4-manifold $(M,g)$ admits three independent compatible complex structures then the metric $g$ is anti-self-dual.
\end{Prop}
By Pontecorvo's classification of ASD bihermitian surfaces \cite{PON} Prop. 3.7, these surfaces are hyperhermitian. Following C.P. Boyer \cite{B} we define {\bf hyperhermitian complex surfaces} as oriented compact conformal 4-manifolds $(M,c,\cal F)$ with a 2-sphere $\cal F$ of compatible complex structures generated by two anti-commuting ones. A hyperhermitian 4-manifold $(M,c,\cal F)$ must be one of the following
\begin{itemize}
\item A flat complex torus
\item A K3 surface with Ricci-flat K\"ahler metric, or
\item A special Hopf surface
\end{itemize}
in particular they all have $b_2(M)=0$.  We refer to \cite{PON} and \cite{B} for details.\\
 Now we focus on the case where first Betti number $b_1(M)$  is odd; for the even case see \cite{AGG}.
\subsection{Numerically anticanonical divisor of a bihermitian surface}
 For the convenience of the reader we recall the results used in the sequel (see \cite{A,AGG,PON}): Denote by 
$$F^g_i(.\ ,\ .)=g(J_i\, .\ , \ . ), \quad i=1,2$$
the K\"ahler forms of $(g,J_1)$ and $(g,J_2)$ respectively, $\t^g_1$, $\t^g_2$ their  Lee forms, i.e.
$$dF^g_i=\t^g_i\wedge F^g_i, \quad i=1,2$$
We furthermore denote by $[J_1,J_2]=J_1J_2-J_2J_1$ the commutator of $J_1$ and $J_2$ and we consider the real $J_i$-anti-invariant $2$-form
$$\Phi^g(.\ ,\ .) = \frac{1}{2}g([J_1,J_2].\ , \ .),$$
and the corresponding complex $(0,2)$-forms
$$\s^g_i(.\ ,\ .) = \Phi^g(.\ ,\ .) + i\Phi^g(J_i.\ ,\ .).$$
Then $\s^g_i$, $i=1,2$ are smooth sections of the anti-canonical bundle $K^{-1}_{J_i}\simeq \bigwedge^{(0,2)}_{J_i}(M)$ of $(M,J_i)$ and $\s^g_i(x)=0$  \iff $\Phi^g(x)=0$ \iff $J_1(x)=\pm J_2(x)$. Therefore the common zero set of $\s_i$ is exactly $\cal D=\cal D^+\cup \cal D^-$, where 
$$\cal D^+=f^{-1}(1)=\{x\in M \mid J_1(x)=J_2(x)\}, $$
$$ \cal D^-=f^{-1}(-1)=\{x\in M \mid J_1(x)=-J_2(x)\}.$$

\begin{Lem}[\cite{A,AGG}]\label{LemmaAGG} Let $(M,c,J_1,J_2)$ be a  bihermitian surface. Then, for any metric $g$ in the conformal class $c$, the $1$-forms $\t^g_1$, $\t^g_2$  and $\s^g_1$ satisfy the following properties:\\
i) If $M$ is compact, then $d(\t^g_1+\t^g_2)=0$,\\
ii) $\bar\partial_{J_1}\s^g_1 = -\frac{1}{2} (\t^g_1 + \t^g_2)^{(0,1)}\otimes \s^g_1$,\\
where $(\ .\ )^{(0,1)}$ denotes the $(0,1)$ part and $\bar\partial_{J_1}$ is the Cauchy-Riemann operator relatively to $J_1$.
\end{Lem}
Proof: i) By \cite{AGG}, $\bigl(d(\t^g_1+\t^g_2)\bigr)_+=0$. Let $\d=-\star d \star$ be the adjoint of $d$. Setting $\f=\t^g_1+\t^g_2$,  $\D d\f=d\d d\f=-d\star d \star d\f=0$, for by assumption $\ast d\f=-d\f$. Hence $\d d\f=0$, and $|d\f |^2=(d\f,d\f)=(\f,\d d\f)=0$.\\
ii) The proof in \cite{AGG} Lemma 3 is local.\hfill $\Box$

\begin{Prop}\label{BNAC} Let $(M,c,J_1,J_2)$ be a compact bihermitian surface. Then there exists a topologically trivial line bundle $L\in H^1(M,\bb R_+^\star)$ such that $\s_i$ is a non-trivial holomorphic section of $K_{J_i}^{-1}\otimes L$, in particular
$$H^0(M,K_{J_i}^{-1}\otimes L)\neq 0,  \quad i=1,2,$$
and $\cal D^+$, $\cal D^-$ are empty or complex curves for both $(M,J_1)$ and $(M,J_2)$.\\
 Moreover, if $(c,J_1,J_2)$ is strongly bihermitian, then $K_{J_i}=L$,
 in particular, $b_2(M)=0$.
 \end{Prop}
Proof: By (\ref{LemmaAGG}), $\t^g_1+\t^g_2$ is closed, hence there is an open covering $(U_j)_{j\in I}$ and $\cal C^\infty$ functions $\Phi_j:U_j\to \bb R$, such that for the local metric $g_j=\exp(\Phi_j)g$, the local Lee forms satisfy $\t^{g_j}_1+\t^{g_j}_2=0$, i.e. $(\t^g_1 +\t^g_2)_{\mid U_j}=-2d\Phi_j$. Setting $c_{jk}=\exp(\Phi_j-\Phi_k)\in \bb R_+^\star$, $\s^j_i=\s^{g_j}_i$, we obtain a topologically trivial line bundle $L=[(c_{jk})]\in H^1(M,\bb R_+^\star)$, and a holomorphic  section $\s_i=(\s^j_i)$ of $K_{J_i}^{-1}\otimes L$. If $(c,J_1,J_2)$ is strongly bihermitian $\cal D=\emptyset$, therefore $\s_i$ is a non-vanishing holomorphic section of $K_{J_i}^{-1}\otimes L$.\hfill $\Box$

\subsection{Bihermitian surfaces with odd first Betti number}
The Kodaira dimension of a compact bihermitian surface with odd first Betti number is $\k=-\infty$ (\cite{A} thm 1), hence the minimal models of $(M,J_1)$ and of $(M,J_2)$ are in class VII$_0$. In this section, we shall complete (and simplify) the classification theorem of V.~Apostolov \cite{A}. We need first

\begin{Lem}\label{BU} Let $S$ be a complex surface and let $\Pi:S'\to S$ be the blowing-up of $x\in S$, and $E=\Pi^{-1}(x)$ its exceptional curve.\\
1) If there exists a flat line bundle $L'$ on $S'$ such that $H^0(S,K_{S'}^{-1}\otimes L')\neq 0$, then for $L=\Pi_\star L'$, $H^0(S,K_S^{-1}\otimes L)\neq 0$.\\
2) If there exists a flat line bundle $L$ on $S$ such that $H^0(S,K_{S}^{-1}\otimes L)\neq 0$, and if $\Pi$ blows-up a point on the effective twisted anticanonical divisor, then for $L'=\Pi^\star L$, $H^0(S',K_{S'}^{-1}\otimes L')\neq 0$.
\end{Lem}
Proof: 1) The coherent sheaf $L=\Pi_\star L'$ is locally trivial since on a simply connected neighbourhood $U$ of the exceptional curve $E$, $L_{\mid U}$ is trivial.  The line bundle $K_S^{-1}\otimes L$ has a section on $S\setminus \Pi(E)$ which extends by Hartogs theorem.\\
2) We have $K_{S'}^{-1}\otimes L'=\pi^\star(K_S^{-1}\otimes L)-E$ hence a section $\s$ of $K_{S}^{-1}\otimes L$ yields a section of $K_{S'}^{-1}\otimes L'$ \iff $x$ belongs to the zero set of $\s$.\hfill $\Box$\\

\begin{Th} Let $(M,c,J_1,J_2)$ be a compact bihermitian surface with odd first Betti number. \\
1) If $(M,c,J_1,J_2)$ is strongly bihermitian (i.e $\cal D=\emptyset$), then the complex surfaces $(M,J_i)$ are minimal and either a Hopf surface covered by a primary one associated to a contraction $F:(\bb C^2,0)\to (\bb C^2,0)$ of the form
$$F(z_1,z_2)=(\a z_1+sz_2^m, a\a^{-1}z_2), $$
$$ {\rm with}\quad a, s\in\bb C, \ 0<|\a|^2\le a< |\a|<1, \ (a^m-\a^{m+1})s=0,$$
or else $(M,J_i)$ are Inoue surfaces $S^+_{N,p,q,r;t}$, $S^-_{N,p,q,r}$.\\
2) If $(M,c,J_1,J_2)$ is not strongly bihermitian, then $\cal D$ has at most two connected components, $(M,J_i)$, $i=1,2$, contain GSS and the minimal models $S_i$ of $(M,J_i)$ are
\begin{itemize}
\item Surfaces with GSS of intermediate type if $\cal D$ has one connected component
\item Hopf surfaces of special type (see {\rm \cite{PON} 2.2}),
Inoue (parabolic) surfaces or Inoue-Hirzebruch surfaces if $\cal D$ has two connected components.
\end{itemize}
Moreover, the blown-up points belong to the NAC divisors.
\end{Th}
Proof:  0) Since the fundamental group of a surface with GSS is isomorphic to $\bb Z$, since $\pi_1( S^{(+)}_{N,p,q,r;t})\neq \bb Z$,   $\pi_1( S^{(-)}_{N,p,q,r})\neq \bb Z$, and  $ S^{(-)}_{N,p,q,r}$ is a quotient of a surface  $ S^{(+)}_{N,p,q,r;t}$ (\cite{I} p 276), $(M,J_1)$ and $(M,J_2)$ must be of the same type.\\
1) If   $(M,c,J_1,J_2)$ is strongly bihermitian, $(M,J_i)$ are minimal using (\ref{BNAC}). Applying the classification of minimal surfaces of  class VII$_0$ with $b_2(S)=0$ \cite{LT}, we have to consider Hopf surfaces and Inoue surfaces. We derive from Proposition (\ref{BNAC}) that $-K_{J_i}\in H^1(M,\bb R_+^\star)$. On one hand the anticanonical line bundle of a surface $S_M$ is not real (\ref{examplesgamma}), hence we may exclude it. On second hand, a finite covering of a bihermitian surface is bihermitian and if a primary Hopf surface is associated to the contraction
$$F(z)=(\alpha_{1}z_{1}+sz_{2}^{m},\alpha_{2}z_{2}), \qquad 0<|\alpha_1|\leq
|\alpha_2|<1,\quad (\alpha_2^m-\alpha_1)s=0,$$
$-K=L^{\a_1\a_2}$, therefore $a=\a_1\a_2 \in \bb R_+^\star$ and it is easy to check that the requested conditions are fulfilled.\\
2) If $(M,c,J_1,J_2)$ is not strongly bihermitian, $(M,J_i)$ admits a non-trivial effective NAC divisor $D_{-K_i}$ whose support is $\cal D$,
and as noticed in \cite{A}, 3.3, it is the same for the minimal model $S_i$ of $(M,J_i)$.  Using theorem (\ref{DNAC}), $S_i$ contains GSS and the support of the NAC divisor has at most two components by (\ref{arbre-IH}). If $\Pi=\Pi_m\cdots\Pi_0:(M,J_i)\to S_i$ is the blowing-down of the exceptional curves, the blown-up points belong to the successive NAC divisors by (\ref{BU}), hence $\cal D$ has at most two connected components.\hfill $\Box$

\begin{Cor}  Let $(M,c,J_1,J_2)$ be a compact  ASD bihermitian surface with odd first Betti number. Then the minimal models of the complex surfaces $(M,J_i)$, $i=1,2$, are 
\begin{itemize}
\item Hopf surfaces of special type (see {\rm \cite{PON} 2.2}),
\item (parabolic) Inoue surfaces or
\item even  Inoue-Hirzebruch surfaces.
\end{itemize}
Moreover, the blown-up points belong to the NAC divisors.
\end{Cor}
Proof: By \cite{PON} 3.11, $(M,J_i)$ has a minimal model $S_i$ in the class VII$_0$ and 
$$H^0(M,-K_{J_i})\neq 0,$$
therefore $S_i$ is not a Inoue surface $S^{(+)}_{N,p,q,r;t}$,   $S^{(-)}_{N,p,q,r}$.
If $b_2(S_i)=0$, $S_i$ is a Hopf surface, and if $b_2(S_i)>0$, $S_i$ contains a GSS by (\ref{DNAC}). The existence of a (non-twisted) global section of $-K$ is equivalent to the existence of a metric $g\in c$ such that $\t^g_1+\t^g_2=0$. In this situation, \cite{AGG} Prop.4 asserts that $\cal D_+$ and $\cal D_-$ are both non-empty. This is possible only when $S_i$ is parabolic Inoue or an even Inoue-Hirzebruch surface.\hfill $\Box$

\begin{flushright}  Centre de Mathématiques et d'Informatique\\
Laboratoire d'Analyse Topologie et Probabilités\\
 Université d'Aix-Marseille 1\\
 39, rue F. Joliot-Curie\\
 13453 Marseille Cedex 13\\
FRANCE\\
\end{flushright}
\end{document}